%% file: main.tex
\documentclass[11pt,twoside]{article}

\usepackage{amsmath, amsthm,amsfonts,amssymb,mathtools,mathdots,scalerel,mathrsfs,stmaryrd,cmll}
\usepackage[utf8]{inputenc}
\usepackage[T1]{fontenc}
\usepackage{mlmodern, tikz-cd, quiver}
\usepackage[many]{tcolorbox}

\usepackage{enumitem}
\usepackage[hidelinks]{hyperref}
\usepackage{tikz, fancyhdr, xparse, xcolor, tocloft}
\usepackage[british]{babel}
\usepackage[a4paper, top=4.5cm, bottom=4.5cm, left=4cm, right=4cm, asymmetric]{geometry}

\usepackage{crossreftools}
\usepackage[textwidth=3cm, textsize=small, colorinlistoftodos]{todonotes}

\usepackage{ahtitle, titlesec}
\usepackage{etoolbox}
\makeatletter
\patchcmd{\ttlh@hang}{\parindent\z@}{\parindent\z@\leavevmode}{}{}
\patchcmd{\ttlh@hang}{\noindent}{}{}{}
\makeatother

\input{macros}

\input{style}

\newcommand\runtitle{equivalences in diagrammatic sets}
\newcommand\runauthor{chanavat and hadzihasanovic}

\title{Equivalences in diagrammatic sets}

\author{Cl\'emence Chanavat and Amar Hadzihasanovic}

\institution{Tallinn University of Technology}

\begin{document}

\thispagestyle{empty}
\maketitle 

\noindent\makebox[\textwidth][r]{%
	\begin{minipage}[t]{.7\textwidth}
\small \emph{Abstract.}
	We show that diagrammatic sets, a topologically sound alternative to polygraphs and strict $\omega$\nbd categories, admit an internal notion of equivalence in the sense of coinductive weak invertibility.
	We prove that equivalences have the expected properties: they include all degenerate cells, are closed under 2-out-of-3, and satisfy an appropriate version of the ``division lemma'', which ensures that enwrapping a diagram with equivalences at all sides is an invertible operation up to higher equivalence.
	On the way to this result, we develop methods, such as an algebraic calculus of natural equivalences, for handling the weak units and unitors which set this framework apart from strict $\omega$\nbd categories.
\end{minipage}}

\vspace{20pt}

\makeaftertitle

\normalsize

\noindent\makebox[\textwidth][c]{%
\begin{minipage}[t]{.55\textwidth}
\setcounter{tocdepth}{1}
\tableofcontents
\end{minipage}}

\input{introduction}
\input{ddiagrams}
\input{equivalences}
\input{natural}

\input{bicategory}
\input{division}

\bibliographystyle{alpha}
\small \bibliography{main.bib}

\end{document}

%% file: macros.tex
\newcommand\eqdef{\coloneqq}
\newcommand\nbd{\nobreakdash-\hspace{0pt}}
\newcommand\idd[1]{\mathrm{id}_{#1}}
\newcommand\invrs[1]{#1^{-1}}
\newcommand\after{\circ}

\newcommand\incl{\hookrightarrow}
\newcommand\iso{\overset{\sim}{\rightarrow}}
\newcommand\isocl[1]{[#1]}
\newcommand\surj{\twoheadrightarrow}
\newcommand\restr[2]{{#1}{\raisebox{0pt}{$|_{#2}$}}}
\newcommand\powerset[1]{\mathscr{P}{#1}}

\newcommand\set[1]{\left\{ {#1} \right\}}
\newcommand\order[2]{#2^{(#1)}}

\newcommand\slice[2]{{#1}/{\raisebox{-2pt}{$#2$}}}

\newcommand\cat[1]{\mathbf{#1}}

\newcommand\fun[1]{\mathsf{#1}}

\newcommand\Pos{\cat{Pos}}

\newcommand\rdcpx{\cat{RDCpx}}
\newcommand{\atom}{{\scalebox{1.3}{\( \odot \)}}}

\newcommand{\Set}{\cat{Set}}

\newcommand{\dgmSet}{\atom\Set}

\newcommand{\pt}{\mathbf{1}}

\DeclareMathOperator{\clos}{cl}
\newcommand\clset[1]{\mathrm{cl}\set{#1}}

\newcommand\gr[2]{#2_{#1}}

\newcommand\bd[2]{\partial_{#1}^{#2}}

\newcommand\faces[2]{\Delta_{#1}^{#2}}
\newcommand\inter[1]{\mathrm{int}\,#1}

\newcommand\cp[1]{\,{\scriptstyle\#}_{#1}\,}
\newcommand\cpsub[1]{\triangleright_{#1}}
\newcommand\subcp[1]{\prescript{}{#1\!}{\triangleleft}\;}
\newcommand\gencp[1]{\,{\scriptstyle\widehat{\#}}_{#1}\,}

\newcommand\subs[3]{#1[#2/#3]}
\newcommand\compos[1]{\langle#1\rangle}

\newcommand\celto{\Rightarrow}
\newcommand\eqvto{\overset{\sim}{\Rightarrow}}

\newcommand\submol{\sqsubseteq}

\newcommand\gray{\otimes}
\newcommand\dual[2]{\fun{D}_{#1}{#2}}

\newcommand{\lcyl}[1]{\mathrm{L}_{#1}}
\newcommand{\rcyl}[1]{\mathrm{R}_{#1}}
\newcommand\pcyl[1]{\ltimes_{#1}}

\DeclareMathOperator{\lydim}{lydim}

\DeclareMathOperator{\pd}{Pd}
\DeclareMathOperator{\rd}{Rd}
\DeclareMathOperator{\eqv}{Eqv}
\DeclareMathOperator{\dgn}{Dgn}
\DeclareMathOperator{\inv}{Inv}
\DeclareMathOperator{\biinv}{BiInv}

\DeclareMathOperator{\ctxeqv}{CtxEqv}
\DeclareMathOperator{\nateqv}{NatEqv}
\DeclareMathOperator{\unitor}{Unitor}
\DeclareMathOperator{\ctxpfw}{CtxPfw}

\newcommand\bic[1]{\mathscr{B}_{#1}}
\newcommand\rev[1]{{#1}^\dag}

\newcommand\fact{\lesssim}

\newcommand\arr{\vec{I}}

\renewcommand{\a}{\alpha}
\renewcommand{\b}{\beta}
\newcommand{\E}{\mathcal{E}}
\newcommand{\I}{\mathcal{I}}
\newcommand{\B}{\mathcal{B}}
\newcommand{\N}{\mathcal{N}}

\newcommand{\T}{\mathcal{T}}

\renewcommand{\th}{\vartheta}

\newcommand{\un}{\varepsilon}
\newcommand{\lun}[1]{\lambda_{#1}}
\newcommand{\run}[1]{\rho_{#1}}
\newcommand{\linv}{\zeta}
\newcommand{\rinv}{\eta}

\newcommand{\R}{\mathscr{R}}

\makeatletter
\newcommand{\oset}[3][0ex]{%
  \mathrel{\mathop{#3}\limits^{
    \vbox to#1{\kern-1.5\ex@
    \hbox{$\scriptstyle#2$}\vss}}}}
\makeatother
\newcommand\qeq{\oset{?}{=}}

%% file: style.tex
\newtheoremstyle{ittheorem}
  {\topsep}   % ABOVESPACE
  {\topsep}   % BELOWSPACE
  {\itshape}  % BODYFONT
  {0pt}       % INDENT (empty value is the same as 0pt)
  {\sffamily \itshape \bfseries} % HEADFONT
  { ---}         % HEADPUNCT
  {5pt plus 1pt minus 1pt} % HEADSPACE
  {}          % CUSTOM-HEAD-SPEC

\newtheoremstyle{itdfn}
  {\topsep}   % ABOVESPACE
  {\topsep}   % BELOWSPACE
  {}  % BODYFONT
  {0pt}       % INDENT (empty value is the same as 0pt)
  {\sffamily \itshape \bfseries} % HEADFONT
  {}         % HEADPUNCT
  {5pt plus 1pt minus 1pt} % HEADSPACE
  {\thmnumber{#2}{\thmnote{\normalfont\ \ %
  {\sffamily(#3)}.}}}          % CUSTOM-HEAD-SPEC

\newtheoremstyle{itrmk}
  {0.5\topsep}   % ABOVESPACE
  {0.5\topsep}   % BELOWSPACE
  {\normalfont}  % BODYFONT
  {0pt}       % INDENT (empty value is the same as 0pt)
  {\sffamily \itshape} % HEADFONT
  { ---}         % HEADPUNCT
  {5pt plus 1pt minus 1pt} % HEADSPACE
  {}          % CUSTOM-HEAD-SPEC
  
\newtheoremstyle{itexm}
  {0.5\topsep}      % ABOVESPACE
  {0.5\topsep}      % BELOWSPACE
  {\normalfont}     % BODYFONT
  {0pt}             % INDENT (empty value is the same as 0pt)
  {\sffamily \itshape \bfseries \color{\mycolor}}        % HEADFONT
  {\\}               % HEADPUNCT
  {5pt plus 1pt minus 1pt} % HEADSPACE
  {\thmname{#1} \thmnumber{#2}{\thmnote{\normalfont\ \ %
  {\sffamily(#3)}.}}}           % CUSTOM-HEAD-SPEC
  
\makeatletter
  \renewcommand\@upn{\textit}
\makeatother

\theoremstyle{ittheorem}
\newtheorem{thm}{Theorem}[section]
\newtheorem*{thm*}{Theorem}
\newtheorem{prop}[thm]{Proposition}
\newtheorem*{prop*}{Proposition}
\newtheorem{cor}[thm]{Corollary}
\newtheorem{lem}[thm]{Lemma}

\theoremstyle{itdfn}
\newtheorem{dfn}[thm]{}
\theoremstyle{itrmk}
\newtheorem{rmk}[thm]{Remark}
\newtheorem{comm}[thm]{Comment}
\newtheorem*{exm}{Example}

\setlist{leftmargin=20pt,itemsep=0pt,parsep=0pt,topsep=1ex}

\titleformat{\section}
 {\large\scshape}{\thesection.}{1em}{}

\titleformat{\subsection}
 {\normalsize\itshape}{\thesubsection.}{1em}{}
\titlespacing*{\subsection}
{0pt}{1.5ex plus 1ex minus .2ex}{1.5ex plus .2ex}

\renewcommand{\cftsecpagefont}{\mdseries}

\makeatletter \renewcommand{\cftsecfillnum}[1]{%
  {\cftsecleader}\nobreak
  \makebox[\@pnumwidth][\cftpnumalign]{\cftsecpagefont \oldstylenums{#1}}\cftsecafterpnum\par
} \makeatother

%% file: introduction.tex
\section*{Introduction}

There is a chicken-and-egg dilemma at the core of higher category theory:
\begin{center}
	\emph{what comes first, higher categories or higher groupoids?}
\end{center}
This is closely connected to the question: \emph{what is an equivalence in a higher category?}, which is, incidentally, the title of a recent survey \cite{ozornova2024equivalence}.
Different models of higher categories reveal different ideological commitments with respect to this question.
Two opposite positions may thus be summarised.
\begin{itemize}
	\item \emph{The homotopist position.}
		The notion of space, or homotopy type, is fundamental.
		It may really be of logical nature, accurately modelling the usage of mathematical equality.
		A higher category, as a structure, has underlying \emph{spaces} of cells.
		An equivalence in a higher category is a cell whose action-by-composition on lower-dimensional cells tracks a homotopy in their underlying space; thus, it is a cell that has an underlying homotopy.
		A higher groupoid is a higher category where every cell tracks an underlying homotopy; thus there is no information in the spaces of higher-dimensional cells that was not already in the space of 0\nbd cells.
		In this sense, higher groupoids \emph{are} really spaces, and have logical priority.
	\item \emph{The computationalist position.}
		Higher groupoids are a special case of higher categories, so the latter must come first. 
		Computation---which is directed and potentially irreversible---is more fundamental than (reversible) homotopy or mathematical equality.
		Whether a cell is or is not an equivalence is about whether it is weakly invertible in a suitable sense, which is an algebraic or computational property.
\end{itemize}
To borrow Girard's terminology \cite{girard2011blind}, the first position entails an \emph{essentialist} view of equivalences---a cell is an equivalence if it \emph{is} a homotopy in the underlying space---while the second entails an \emph{existentialist} view: a cell is an equivalence if it \emph{behaves} like one.

Subscription to one or the other position matches quite neatly the divide between \emph{non-algebraic} and \emph{algebraic} models of higher categories.
In most well-established non-algebraic models, a higher category is equipped with an underlying space which is extra data not definable in the algebraic language of cells, units, and composition.
In the complete Segal $\Theta_n$-space model \cite{rezk2010cartesian}, among others, this is achieved by working directly with space-valued presheaves.
In the complicial model \cite{verity2008complicial}, this is achieved by equipping the underlying simplicial set of a higher category with a ``marked'' subset of equivalences.
Since these models typically generalise a model of higher groupoids which is already known to be sound for classical homotopy theory, the \emph{homotopy hypothesis}---that higher groupoids model all classical homotopy types---is usually an established theorem.

Algebraic models, on the other hand, by their nature focus on notions of equivalence that are definable in the algebra of units and composition.
In particular, much attention has been given to the coinductive notion of \emph{pseudoinvertible cell} \cite{cheng2007omega}, also known as ``weakly invertible cell'' or ``$\omega$\nbd equivalence''.
Originally confined to strict $\omega$\nbd categories, where it plays a crucial role in the definition of the folk model structure \cite{lafont2010folk}, this notion has recently been studied also in the context of weak algebraic models \cite{rice2020coinductive, fujii2024weakly, benjamin2024invertible}.
Beyond philosophy, there are technical reasons why definable, algebraic equivalences are convenient.
Most of these are linked to the existence of the \emph{walking $\omega$\nbd equivalence}, a classifying object for equivalences, and in particular its \emph{coherent} or contractible version, for which an explicit cellular model was exhibited in \cite{hadzihasanovic2025model}.
A description of this object provides an explicit, computable model for the localisation of a higher category at a cell or set of cells.

In recent years, a growing network of equivalences between different non-algebraic models \cite{ozornova2023quillen, doherty2023equivalence, loubaton2023theory} and between different algebraic models \cite{ara2010groupoides, bourke2020iterated, benjamin2024invertible} has been produced, but notably no equivalence between a non-algebraic and an algebraic model.
More in general, the homotopy hypothesis has not been proven for any of the algebraic models.
(One partial exception is \cite{henry2018regular}, but this is only partially algebraic---composition is algebraic, but units are not---and only a model of higher groupoids.)
Indeed, even producing model structures on categories of algebraic higher categories that can plausibly model the homotopy theory of $(\infty, n)$-categories has proven challenging; the works studying equivalences in and between algebraic weak $n$\nbd categories can be seen as steps in this direction.

Following the second-named author's lead in the unpublished \cite{hadzihasanovic2020diagrammatic}, we have started in \cite{chanavat2024diagrammatic} a programme to develop a model of higher categories that could serve as a ``bridge'' between the non-algebraic and algebraic models, centred on the notion of \emph{diagrammatic set}.
Like cubical and simplicial models, this is based on presheaves on a shape category, whose objects are combinatorial models of topological closed balls.
This allowed us, in our first article, to use Cisinski's methods \cite{cisinski2006prefaisceaux} in order to prove the homotopy hypothesis for our model of higher groupoids.

On the other hand, the rich combinatorics of diagram shapes available in diagrammatic sets allows for modes of reasoning very close to those available in strict $n$\nbd categories.
Indeed, the original motivation for diagrammatic sets was to provide a ``topologically sound'' alternative to strict $n$\nbd categories and polygraphs \cite{ara2023polygraphs} for the purposes of higher-dimensional diagram rewriting.
The main difference between diagrammatic sets and polygraphs is that the input and output of a cell in a diagrammatic set must be \emph{round} diagrams, that is, diagrams that are ``shaped'' like topological balls of the appropriate dimension.
This is the same restriction that appears in Henry's regular polygraphs \cite{henry2018regular}, but unlike regular polygraphs, diagrammatic sets have a rich algebra of \emph{weak units} and other degenerate cells which can be used to ``pad'' diagrams until they are round:

\[\begin{tikzcd}[column sep=small]
	&&&&&&& {{\scriptstyle x}\;\bullet} \\
	\bullet & {{\scriptstyle x}\;\bullet} && \bullet && \bullet & {{\scriptstyle x}\;\bullet} && \bullet \\
	{\text{not round}} && \bullet &&& {\text{round}} && \bullet
	\arrow["g", curve={height=-6pt}, from=1-8, to=2-9]
	\arrow["f", from=2-1, to=2-2]
	\arrow[""{name=0, anchor=center, inner sep=0}, "g", from=2-2, to=2-4]
	\arrow[curve={height=6pt}, from=2-2, to=3-3]
	\arrow[""{name=1, anchor=center, inner sep=0}, "f", curve={height=-12pt}, from=2-6, to=1-8]
	\arrow["f", from=2-6, to=2-7]
	\arrow["x"{description}, from=2-7, to=1-8]
	\arrow[""{name=2, anchor=center, inner sep=0}, "g"{description}, from=2-7, to=2-9]
	\arrow[curve={height=6pt}, from=2-7, to=3-8]
	\arrow[curve={height=6pt}, from=3-3, to=2-4]
	\arrow[curve={height=6pt}, from=3-8, to=2-9]
	\arrow["f"', shorten >=3pt, Rightarrow, from=2-7, to=1]
	\arrow["g"', shorten <=3pt, Rightarrow, from=2, to=1-8]
	\arrow[shorten >=3pt, Rightarrow, from=3-3, to=0]
	\arrow[shorten >=3pt, Rightarrow, from=3-8, to=2]
\end{tikzcd}\]
In this article, we show that, like the algebraic models, diagrammatic sets support a coinductive definition of equivalence---in fact, several equivalent definitions, one of which is in terms of ``pseudoinvertibility''---with the good properties that one expects.
This will be a crucial step towards the definition of model structures for $(\infty, n)$\nbd categories on the category of diagrammatic sets, that will mix attributes of the folk model structure on strict $n$\nbd categories and of presheaf model structures \`a la Cisinski.

The key result about coinductive equivalences used in the construction of the folk model structure is the \emph{division lemma} \cite[Lemma 20.1.10]{ara2023polygraphs}.
This states, roughly, that the action by $k$\nbd composition of a $(k+1)$\nbd equivalence induces a bijection of sets of $n$\nbd cells up to $(n+1)$\nbd equivalence for each $n > k$.
When trying to import this result into our setting, one sees that there is an immediate combinatorial obstacle: pasting at the $k$\nbd boundary only sends round $n$\nbd dimensional diagrams to round $n$\nbd dimensional diagrams when $k = n-1$.
For example, when $k = 0$ and $n = 2$, pasting at the input boundary looks like
\[\begin{tikzcd}
	\bullet & \bullet & \mapsto & \bullet & \bullet & \bullet
	\arrow[""{name=0, anchor=center, inner sep=0}, curve={height=-12pt}, from=1-1, to=1-2]
	\arrow[""{name=1, anchor=center, inner sep=0}, curve={height=12pt}, from=1-1, to=1-2]
	\arrow[from=1-4, to=1-5]
	\arrow[""{name=2, anchor=center, inner sep=0}, curve={height=12pt}, from=1-5, to=1-6]
	\arrow[""{name=3, anchor=center, inner sep=0}, curve={height=-12pt}, from=1-5, to=1-6]
	\arrow[shorten <=3pt, shorten >=3pt, Rightarrow, from=1, to=0]
	\arrow[shorten <=3pt, shorten >=3pt, Rightarrow, from=2, to=3]
\end{tikzcd}\]
and the resulting non-round diagram cannot appear as the input or output of an equivalence.
To obtain a well-formed statement, one needs to replace the action-by-composition of a single diagram with the action of a \emph{round context}, that is, a ``round diagram with a round hole'', whose action may look like
\[\begin{tikzcd}
	&&&&& \bullet \\
	\bullet & \bullet & \mapsto & \bullet & \bullet & \bullet & \bullet \\
	&&&& \bullet
	\arrow[curve={height=-6pt}, from=1-6, to=2-7]
	\arrow[""{name=0, anchor=center, inner sep=0}, curve={height=-12pt}, from=2-1, to=2-2]
	\arrow[""{name=1, anchor=center, inner sep=0}, curve={height=12pt}, from=2-1, to=2-2]
	\arrow[""{name=2, anchor=center, inner sep=0}, curve={height=-12pt}, from=2-4, to=1-6]
	\arrow[from=2-4, to=2-5]
	\arrow[curve={height=6pt}, from=2-4, to=3-5]
	\arrow[curve={height=-6pt}, from=2-5, to=1-6]
	\arrow[""{name=3, anchor=center, inner sep=0}, curve={height=12pt}, from=2-5, to=2-6]
	\arrow[""{name=4, anchor=center, inner sep=0}, curve={height=-12pt}, from=2-5, to=2-6]
	\arrow[Rightarrow, from=2-6, to=1-6]
	\arrow[from=2-6, to=2-7]
	\arrow[Rightarrow, from=3-5, to=2-5]
	\arrow[curve={height=6pt}, from=3-5, to=2-6]
	\arrow[shorten <=3pt, shorten >=3pt, Rightarrow, from=1, to=0]
	\arrow[shorten >=3pt, Rightarrow, from=2-5, to=2]
	\arrow[shorten <=3pt, shorten >=3pt, Rightarrow, from=3, to=4]
\end{tikzcd}\]
so that round diagrams are mapped to round diagrams.
Our main result is a proof of the division lemma for \emph{weakly invertible} round contexts, which are ``built out of equivalences'' in an appropriate sense.
A special case are pastings of equivalences at lower-dimensional boundaries which are ``rounded'' with the use of weak units, as sketched above.
Note that, in the framework of strict $n$\nbd categories, one can recover the original form of the division lemma as a special case, by rounding with \emph{strict} rather than weak units.

At the outset, we tried to develop a minimum of methods for handling weak units that would allow us to reproduce the original proof of the division lemma for strict $n$\nbd categories.
This proof used a number of ``tricks'' and explicit diagrammatic calculations.
In the process, we have developed a number of abstractions---in particular, an algebraic calculus of \emph{natural equivalences} of round contexts, where naturality is itself coinductively defined---which have resulted in an new, uniform, higher-level proof.
We believe that these methods may also be of interest in the theory of strict $n$\nbd categories.

%%%%%%%%%%%

\subsection*{Background on diagrammatic sets}

The theory of diagrammatic sets is built upon the combinatorics of \emph{molecules} and of \emph{regular directed complexes}.
Here we give a brief overview, and refer to the book \cite{hadzihasanovic2024combinatorics} for details.

It is a classical result of combinatorial topology that \emph{regular cell complexes} can be reconstructed up to cellular homeomorphism from their \emph{face poset}, which is graded by dimension.
Intuitively, the shape of a higher-categorical diagram is a \emph{directed cell complex}, that is, a cell complex in which the $n$\nbd dimensional faces of an $(n+1)$\nbd dimensional cell $x$ are partitioned into an \emph{input} half $\faces{}{-}x$ and an \emph{output} half $\faces{}{+}x$.
The notion of \emph{regular directed complex} attempts to capture a class of higher-categorical diagram shapes that are fully described by their \emph{oriented face poset}, which records this extra information.

Formally, one starts with a category of \emph{oriented graded posets}, whose objects are graded posets $P = \bigcup_{n \in \mathbb{N}} \gr{n}{P}$ together with a bipartition $\faces{}{}x = \faces{}{+}x + \faces{}{-}x$ of the set of faces (covered elements) of each $x \in P$, together with morphisms $f\colon P \to Q$ that induce a bijection between $\faces{}{\a}x$ and $\faces{}{\a}f(x)$ for each $\a \in \set{+, -}$ and $x \in P$.
Oriented graded posets can be equipped with the Alexandrov topology, in which closed sets are those that contain the lower set of each of their elements; all morphisms are both closed and continuous with respect to this topology.
We write $\clos$ for its closure operator.

Given a closed subset $U \subseteq P$, let $\dim{U}$ be $-1$ if $U = \varnothing$, the maximal dimension of an element of $U$ if one exists, and $\infty$ otherwise.
For each $n \geq -1$, $U$ admits a notion of \emph{input $n$\nbd boundary} $\bd{n}{-}U$ and \emph{output $n$\nbd boundary} $\bd{n}{+}U$, both of which are closed subsets of dimension $\leq n$.
We usually omit $n$ when it is equal to $\dim{U} - 1$, and let $\bd{n}{}U \eqdef \bd{n}{+}U \cup \bd{n}{-}U$.
We say that $U$ is \emph{globular} if $\bd{k}{\a}\bd{n}{\b}U = \bd{k}{\a}U$ for all $n \in \mathbb{N}$, $k < n$ and $\a, \b \in \set{+, -}$, and \emph{round} if it also satisfies $\bd{n}{+}U \cap \bd{n}{-}U = \bd{n-1}{}U$ for all $n < \dim{U}$.

\emph{Molecules} are a subclass of oriented graded posets aiming to capture the shapes of regular \emph{pasting diagrams}, that is, those diagrams that are composable in the algebra of strict $\omega$\nbd categories.
They are generated inductively by the following three clauses.
First of all, the oriented graded poset $\pt$ with a single element and trivial orientation is a molecule: this is the \emph{point}, or the shape of a 0\nbd cell.
Next, if $U$ and $V$ are molecules, and there exists an isomorphism $\bd{k}{+}U \iso \bd{k}{-}V$, then the oriented graded poset $U \cp{k} V$ obtained as the pushout of $\bd{k}{+}U \incl U$ and $\bd{k}{+}U \iso \bd{k}{-}V \incl V$ is a molecule, the \emph{pasting of $U$ and $V$ at the $k$\nbd boundary}.
Finally, if $U$ and $V$ are \emph{round} molecules of the same dimension $n$, and we have isomorphisms $\bd{}{\a}U \iso \bd{}{\a}V$ for all $\a \in \set{+, -}$, which determine an isomorphism $\bd{}{}U \iso \bd{}{}V$, then the oriented graded poset $U \celto V$ obtained by first taking the pushout of $\bd{}{}U \incl U$ and $\bd{}{}U \iso \bd{}{}V \incl V$, then adding a new greatest element $\top$ with $\faces{}{-}\top = \gr{n}{U}$ and $\faces{}{+}\top = \gr{n}{V}$, is a molecule.
This final operation should be thought of as the construction of a closed $(n+1)$\nbd ball by first gluing two closed $n$\nbd balls along their boundary to obtain an $n$\nbd sphere, then filling it with an open $(n+1)$\nbd ball.

Molecules have a number of remarkable properties.
They are rigid, in that they do not have any non-trivial automorphisms; this also implies that $U \cp{k} V$ and $U \celto V$ are independent of the boundary isomorphisms in their definition.
They are all globular, and, in fact, their isomorphism classes form a strict $\omega$\nbd category with pasting at the $k$\nbd boundary as $k$\nbd composition.
Not all molecules are round: roundness coincides, for molecules, with the topological property of having an order complex homeomorphic to a closed ball.

Most molecules admit many non-trivial pasting decompositions, but some order can be found via the notion of \emph{$k$\nbd layering}, which is a decomposition $U = \order{1}{U} \cp{k} \ldots \cp{k} \order{m}{U}$ with the property that each $\order{i}{U}$ contains a single maximal element $\order{i}{x}$ of dimension $> k$.
To each molecule $U$, one can associate an integer $-1 \leq \lydim{U} < \dim{U}$, the \emph{layering dimension}, such that $U$ is guaranteed to admit a $k$\nbd layering for each $\lydim{U} \leq k < \dim{U}$; moreover, each factor in a $k$\nbd layering has strictly smaller layering dimension, which allows us to prove statements about molecules by \emph{induction on layering dimension}.

The \emph{submolecule inclusions} are the class of morphisms of molecules generated by the inclusions into pastings $U, V \incl U \cp{k} V$ and closed under isomorphisms and composition.
Molecules are closed under the following generalisation of pasting, that we call \emph{pasting at a submolecule}: given a submolecule inclusion $\iota\colon \bd{k}{+}U \incl \bd{k}{-}V$, the pushout $U \cpsub{k,\iota} V$ of $\bd{k}{+}U \incl U$ and $\bd{k}{+}U \stackrel{\iota}{\incl} \bd{k}{-}V \incl V$ is a molecule.
Dually, given $\iota\colon \bd{k}{-}V \incl \bd{k}{+}U$, the pushout $U \subcp{k,\iota} V$ of $\bd{k}{-}V \stackrel{\iota}{\incl} \bd{k}{+}U \incl U$ and $\bd{k}{-}V \incl V$ is a molecule.

If a closed subset inclusion $V \subseteq U$ is a submolecule inclusion, we write $V \submol U$ and say that $V$ is a submolecule of $U$.
In particular, if $\dim{V} = \dim{U}$ and $V$ is round, we say that $V$ is a \emph{rewritable submolecule}.
Rewritable submolecules have the property that one can remove the interior $\inter{V} \eqdef V \setminus \bd{}{}V$ of $V$ from $U$ and replace it with the interior of another round molecule $W$ with the same boundaries, to obtain a new molecule $\subs{U}{W}{V}$.

An \emph{atom} is a molecule with a greatest element; this is either the point, or isomorphic to $U \celto V$ for some round molecules $U$, $V$.
A \emph{regular directed complex} is an oriented graded poset $P$ with the property that the lower set $\clset{x}$ of each $x \in P$ is an atom; all molecules are regular directed complexes.
A \emph{map} $f\colon P \to Q$ of regular directed complexes is an order-preserving function of the underlying posets with the following property: for all $x \in P$, $n \in \mathbb{N}$, and $\a \in \set{+, -}$, we have $f(\bd{n}{\a}x) = \bd{n}{\a}f(x)$, and, furthermore, the restriction $\restr{f}{\bd{n}{a}x}$ is \emph{final} onto its image; in this context, this means that, for all \( y, y' \in \bd{n}{\a}x \), if \( f(y) = f(y') \), then there exists a zig-zag \( y \leq y_1 \geq \ldots \leq y_m \geq y' \) in \( \bd{n}{\a}x \) such that \( f(y) \leq f(y_i) \) for all \( y \in \set{1, \ldots, m} \).
Maps are closed and dimension-non-increasing; a map which preserves the dimension of all elements is, equivalently, a morphism of the underlying oriented graded posets.

While seemingly technical, maps are characterised among order-preserving maps of the underlying posets by the property that they admit a natural interpretation as strict functors of strict $\omega$\nbd categories.
In \cite{chanavat2024diagrammatic}, we further restricted our attention to \emph{cartesian maps}: maps that are, additionally, Grothendieck fibrations of the underlying posets.
We let $\rdcpx$ denote the category of regular directed complexes and cartesian maps, and $\atom$ be a skeleton of its full subcategory on the atoms.
Both categories admit a number of interesting functorial operations: for the purposes of this article, we mention \emph{Gray products} $P, Q \mapsto P \gray Q$, which are a directed version of cartesian products and determine a semicartesian monoidal structure on both categories, and \emph{duals} $P \mapsto \dual{n}{P}$, which leave the underlying poset unchanged but reverse the direction of $n$\nbd dimensional faces.

A \emph{diagrammatic set} is a presheaf on $\atom$; with their morphisms of presheaves, diagrammatic sets form a category $\dgmSet$.
In \cite{chanavat2024diagrammatic}, we proved that $\atom$ with its natural grading is an Eilenberg--Zilber category; this implies that every diagrammatic set is a ``cell complex'' built by attaching atoms along their boundary in successive dimensions.
We also proved that $\rdcpx$ can be identified with the full subcategory of $\dgmSet$ on the ``regular cell complexes'', whose attaching maps are monomorphisms.
Throughout this article, we will often identify a regular directed complex with its representation in $\dgmSet$.

%%%%%%%%%%% 

\subsection*{Structure of the article}

We start in Section \ref{sec:diagrams} by setting up some terminology and notation relative to diagrams in a diagrammatic set, which are morphisms whose domain is a regular directed complex.
In particular, we define \emph{pasting diagrams} and \emph{round diagrams} to be the diagrams whose domain is a molecule and a round molecule, respectively, and a \emph{subdiagram} of a pasting diagram to be its restriction along a submolecule inclusion.
Then, we focus on \emph{degenerate diagrams}, which are those that factor through a cartesian map that strictly decreases dimension.
We prove that each pasting diagram admits certain useful degenerate diagrams living on top of it: \emph{units}, which ``raise'' a diagram to the next dimension, as well as \emph{unitors}, which introduce or eliminate units at a subdiagram of its boundary.
Moreover, degenerate diagrams can be \emph{reversed}, and degenerate pasting diagrams admit further degenerate pasting diagrams, the \emph{invertors}, which exhibit their reverse as an ``inverse up to a higher-dimensional degenerate diagram''.

In Section \ref{sec:equivalences}, we define a coinductive subclass of the round diagrams in a diagrammatic set, whose members we call \emph{equivalences}.
We first give a definition in terms of the existence of ``lax solutions''---that is, solutions exhibited by a round diagram one dimension higher---to certain equations of round diagrams.
Later, we prove that it is equivalent to one given in terms of weak invertibility and one given in terms of ``bi-invertibility'', that is, existence of a separate left and right weak inverse (Theorem \ref{thm:eqv_equals_inv_equals_biinv}).
These three definitions can be seen as weakenings of three equivalent characterisations of isomorphisms $e\colon a \iso b$ in a category, respectively:
\begin{enumerate}
	\item morphisms $e\colon a \to b$ such that, for all morphisms $f$ with codomain $b$ and $g$ with domain $a$, the equations $e \after x \qeq f$ and $x \after e \qeq g$ admit a solution;
	\item morphisms $e\colon a \to b$ such that there exists a morphism $e^*\colon b \to a$ satisfying $e^* \after e = \idd{a}$ and $e \after e^* = \idd{b}$;
	\item morphisms $e\colon a \to b$ such that there exist morphisms $e^L, e^R\colon b \to a$ satisfying $e^L \after e = \idd{a}$ and $e \after e^R = \idd{b}$.
\end{enumerate}
We prove that all degenerate round diagrams are equivalences, and that equivalences are closed under a suitable form of ``2-out-of-3'' (Theorem \ref{thm:equivalences_2_out_of_3}).

Section \ref{sec:natural} begins with the definitions of contexts---in particular, round contexts and weakly invertible contexts---for pasting diagrams in a diagrammatic set.
After proving certain factorisation results for contexts, analogous to the existence of layerings for pasting diagrams, we define the key notion of \emph{natural equivalence} of round contexts.
A natural equivalence is a family of equivalences, indexed by round diagrams in the domain of a context, which satisfy a naturality condition up to higher natural equivalence.
The rest of the section is devoted to the proof that natural equivalences satisfy closure properties analogous to those of natural isomorphisms of functors (Theorem \ref{thm:natural_saturation}): they compose, they can be ``whiskered'' with round contexts on the left and on the right, and they can be inverted.
Moreover, the families of units and unitors defined in Section \ref{sec:natural} are natural in their parameters, and so is the application of an $(n+1)$\nbd dimensional equivalence to an appropriate subdiagram of an $n$\nbd dimensional context.
On the whole, these results determine a kind of algebraic calculus which can be used to construct natural equivalences of round contexts from basic building blocks which, crucially, include unitors; this calculus will be our main tool for manipulating weak units.

Section \ref{sec:bicategory} is a sort of interlude: we show that, if we restrict to the ``slice'' of a diagrammatic set consisting of round diagrams in dimensions $n, n+1, n+2$, and quotient the latter under $(n+3)$\nbd dimensional equivalence, we can give the resulting 2\nbd graph the structure of a bicategory (Proposition \ref{prop:bicategory_well_defined}).
The main interest of this result is that it allows us to import general results about bicategories into the theory of diagrammatic sets---for example, the fact that every equivalence can be promoted to an adjoint equivalence (Proposition \ref{prop:adjointification})---and justifies the use of the calculus of string diagrams for proofs that only involve three consecutive dimensions.
None of this is used later in the article, and is only included because of thematic resonance.

In Section \ref{sec:division}, we formalise the idea of ``rounding'' a context by padding it with weak units, and show how this operation interacts with natural equivalence of contexts, the identity context, and composition of contexts.
Finally, we assemble together all this machinery for our main result, which retroactively justifies a piece of terminology: a weakly invertible round context, is, indeed, invertible up to natural equivalence (Theorem \ref{thm:invertible_context_are_invertible}).
The division lemma (Lemma \ref{lem:division_lemma}) appears as an immediate corollary.

%%%%%%%%%%%

\subsection*{Acknowledgements}

The second-named author was supported by Estonian Research Council grant PSG764.
We thank Yuki Maehara for helpful discussions about diagram paddings and naturality.

%% file: ddiagrams.tex
\section{Diagrams in diagrammatic sets} \label{sec:diagrams}

\subsection{Basic definitions}

\begin{dfn}[Diagram in a diagrammatic set]
	Let $U$ be a regular directed complex and $X$ a diagrammatic set.
	A \emph{diagram of shape $U$ in $X$} is a morphism $u\colon U \to X$.
	A diagram is a \emph{pasting diagram} if $U$ is a molecule, a \emph{round diagram} if $U$ is a round molecule, and a \emph{cell} if $U$ is an atom.
	We write $\dim{u} \eqdef \dim{U}$.
\end{dfn}

\begin{rmk}
	By the Yoneda lemma, a cell in $X$ is the same as an element of $X$ as a presheaf.
\end{rmk}

\begin{rmk}
	Every cell is a round diagram and every round diagram is a pasting diagram.
	Since isomorphisms of molecules are unique when they exist, we can safely identify pasting diagrams that are isomorphic in the slice of $\dgmSet$ over $X$.
	Equivalently, we may assume to have fixed a skeleton of the full subcategory of $\rdcpx$ on the molecules; see \cite{hadzihasanovic2023higher} for an explicit encoding of isomorphism classes of molecules.
\end{rmk}

\noindent Recall that an \emph{$\omega$\nbd graph}, or globular set, is a graded set $G = \sum_{n \in \mathbb{N}} \gr{n}{G}$ together with boundary functions $\bd{}{-}, \bd{}{+}\colon \gr{n+1}{G} \to \gr{n}{G}$ for each $n \in \mathbb{N}$, satisfying $\bd{}{\a}\bd{}{+} = \bd{}{\a}\bd{}{-}$ for all $\a \in \set{+, -}$.
Given $n > 0$ and $a \in G_n$, we write $a\colon a^- \celto a^+$ to express that $\bd{}{-}a = a^-$ and $\bd{}{+}a = a^+$, and say that $a$ is of \emph{type $a^- \celto a^+$}.
We say that $a, b \in G_n$ are \emph{parallel} if $n = 0$, or $n > 0$ and $a, b$ have the same type.

More in general, for each $k \in \mathbb{N}$,  we consider \emph{$\omega$\nbd graphs in degree $\geq k$}; these are no different except for a shift in the indexing.
Given an $\omega$\nbd graph $G$, $k \leq n \in \mathbb{N}$, and $\a \in \set{+, -}$, we let $\bd{k}{\a}\colon \gr{n}{G} \to \gr{k}{G}$ be defined recursively by $\bd{k}{\a} \eqdef \idd{\gr{k}{G}}$ if $n = k$ and $\bd{k}{\a} \eqdef \bd{}{\a}\bd{k+1}{\a}$ if $n > k$.
Given parallel $a, b \in \gr{k}{G}$, 
\[
	G(a, b) \eqdef \set{c \in \sum_{n > k} \gr{n}{G} \mid \bd{k}{-}c = a, \bd{k}{+}c = b}
\]
admits a structure of $\omega$\nbd graph in degree $> k$ with the same grading and boundary functions as $G$.

\begin{dfn}[The $\omega$-graph of pasting diagrams]
	Let $u\colon U \to X$ be a pasting diagram in a diagrammatic set, $n \in \mathbb{N}$, and $\a \in \set{ +, - }$.
	We let $\bd{n}{\a}u \eqdef \restr{u}{\bd{n}{\a}{U}}\colon \bd{n}{\a}U \to X$.
	We may omit the index $n$ when $n = \dim{u} - 1$.

	We let $\pd X$ denote the set of pasting diagrams in $X$ and $\rd X \subset \pd X$ its subset of round diagrams.
	The set $\pd X$ is graded by dimension; given a subset $A$ of $\pd X$ and $n \in \mathbb{N}$, we let $\gr{n}{A} \eqdef \set{u \in A \mid \dim u = n}$.
	Then, $\pd X$ admits the structure of an $\omega$\nbd graph with the functions $\bd{}{-}, \bd{}{+}\colon \gr{n+1}{\pd X} \to \gr{n}{\pd X}$ for each $n \in \mathbb{N}$.
	These restrict along the inclusions $\gr{n}{\rd X} \subseteq \gr{n}{\pd X}$, making $\rd X$ an $\omega$\nbd subgraph of $\pd X$.
\end{dfn}

\noindent	Combining these notations, if \( u, v \) is a pair of \( k \)\nbd dimensional parallel pasting diagrams in a diagrammatic set \( X \) and \( n > k \), then \( \gr{n}{\pd X(u, v)} \) is the set of pasting diagrams \( w \in \pd X \) such that \( \dim w = n \), \( \bd{k}{-} w = u \), and \( \bd{k}{+} w = v \).
	In particular, \( \gr{k + 1}{\pd X(u, v)} \) is the set of pasting diagrams \( h \colon u \celto v \).
	Notice that, if \( u \) and \( v \) are of shape \( U \) and \( V \), respectively, then \( h \colon u \celto v \) need \emph{not} be of shape \( U \celto V \); the shape of \( h \) is not necessarily an atom.

\begin{dfn}[Subdiagram]
	Let $u\colon U \to X$ be a pasting diagram.
	A \emph{subdiagram of $u$} is a pair of 
	\begin{enumerate}
		\item a pasting diagram $v\colon V \to X$, and 
		\item a submolecule inclusion $\iota\colon V \incl U$
	\end{enumerate}
 	such that $v = u \after \iota$.
	A subdiagram is \emph{rewritable} when $\iota$ is a rewritable submolecule inclusion.
	We write $\iota\colon v \submol u$ for the data of a subdiagram of $u$.
\end{dfn}

\noindent We will simply write $v \submol u$ when $\iota$ is irrelevant or evident from the context.

\begin{dfn}[Pasting of pasting diagrams]
	Let $u\colon U \to X$ and $v\colon V \to X$ be pasting diagrams such that $\bd{k}{+}u = \bd{k}{-}v$.
	We let $u \cp{k} v\colon U \cp{k} V \to X$ be the pasting diagram determined by the universal property of the pasting $U \cp{k} V$.

	More generally, suppose we have a subdiagram $\iota\colon \bd{k}{+}u \submol \bd{k}{-}v$.
	We let $u \cpsub{k,\iota} v\colon U \cpsub{k,\iota} V \to X$ be the pasting diagram determined by the universal property of $U \cpsub{k,\iota} V$ as a pasting of $U$ at a submolecule of $\bd{k}{-} V$.
	Dually, if $\iota\colon \bd{k}{-}v \submol \bd{k}{+}u$, we let $u \subcp{k,\iota} v$ be the universally determined pasting diagram of shape $U \subcp{k,\iota} V$.
\end{dfn}

\noindent We may omit the index $k$ when it is equal to $\min \set{\dim{u}, \dim{v}} - 1$, and omit $\iota$ when it is irrelevant or evident from the context. 
Thus, \( u \cp{} v \) stands for \( u \cp{k} v \), with \( k \eqdef \min \set{\dim{u}, \dim{v}} - 1 \). 

\begin{rmk}
	When $\iota$ is an isomorphism, we have $u \cpsub{k,\iota} v = u \subcp{k,\iota} v = u \cp{k} v$.
\end{rmk}

\begin{rmk}
	There are evident subdiagrams $u, v \submol u \cpsub{k,\iota} v$ and $u, v \submol u \subcp{k,\iota} v$ whenever the pastings are defined.
\end{rmk}

\begin{rmk} \label{rmk:strict_omegacats}
	It follows from the results of \cite[Chapter 5]{hadzihasanovic2024combinatorics} that pasting satisfies all the axioms of composition in strict $\omega$\nbd categories.
	In particular, pastings of the form $u \cp{} v$ suffice to generate all pastings of the form $u \cp{k} v$, as well as pastings at a subdiagram $u \cpsub{k, \iota} v$, for all $k \in \mathbb{N}$.
\end{rmk}

\begin{dfn}[Substitution at a rewritable subdiagram]
	Let $u\colon U \to X$ be a pasting diagram, let $\iota\colon v \submol u$ be a rewritable subdiagram of shape $V$, and let $w$ be a round diagram of shape $W$, parallel to $v$.
	The \emph{substitution of $w$ for $\iota\colon v \submol u$} is the unique pasting diagram $\subs{u}{w}{\iota(v)}$ of shape $\subs{U}{W}{\iota(V)}$ which restricts to $w$ along $W \incl \subs{U}{W}{\iota(V)}$ and to $\restr{u}{U \setminus \inter{\iota(V)}}$ along $U \setminus \inter{\iota(V)} \incl \subs{U}{W}{\iota(V)}$.
\end{dfn}

\begin{rmk}
	By \cite[Lemma 7.1.9]{hadzihasanovic2024combinatorics}, if $\dim{u} = \dim{v}$ and $v$ is round, then whenever $u \cpsub{\iota} v$ is defined, it is round.
	Moreover, if $u$ is also round, by \cite[Lemma 7.1.10]{hadzihasanovic2024combinatorics} we have $\bd{}{-}(u \cpsub{\iota} v) = \subs{\bd{}{-}v}{\bd{}{-}u}{\iota(\bd{}{+}u)}$.
	Dual facts hold for $u \subcp{\iota} v$, when defined.
\end{rmk}

%%%%%%%%%%%%

\subsection{Degenerate diagrams}

\begin{dfn}[Degenerate diagram]
	Let $u\colon U \to X$ be a diagram in a diagrammatic set.
	We say that $u$ is \emph{degenerate} if there exists a diagram $v\colon V \to X$ and a surjective cartesian map of regular directed complexes $p\colon U \surj V$ such that $u = v \after p$ and $\dim{v} < \dim{u}$.
\end{dfn}

\begin{dfn}[Reverse of a degenerate diagram]
	Let $u\colon U \to X$ be a degenerate diagram in a diagrammatic set, equal to $v \after p$ for some diagram $v\colon V \to X$ and surjective map $p\colon U \surj V$ with $n \eqdef \dim{u} > \dim{v}$.
	The \emph{reverse of $u$} is the degenerate diagram $\rev{u} \eqdef v \after \dual{n}{p}$ of shape $\dual{n}{U}$.
\end{dfn}

\begin{rmk}
	For each $x \in \gr{n}{U}$, by \cite[Lemma 2.3]{chanavat2024diagrammatic} there is a canonical Eilenberg--Zilber factorisation $(p_x, v_x)$ of the cell $\restr{u}{\clset{x}}$, such that $\restr{\rev{u}}{\clset{\dual{n}{x}}} = v_x \after \dual{n}{p_x}$.
	Since $x$ is arbitrary in $\gr{n}{U}$, and $\dual{n}{}$ acts trivially on lower-dimensional elements, the reverse of $u$ is independent of the choice of factorisation.
\end{rmk}

\begin{rmk}
	If $u$ is a degenerate pasting diagram of type $u^- \celto u^+$, then $\rev{u}$ is of type $u^+ \celto u^-$.
\end{rmk}

\noindent The aim of this section is to introduce some useful families of degenerate pasting diagrams that always exist in a diagrammatic set.

\begin{dfn}[Partial cylinder]
Let $\arr \eqdef \pt \celto \pt$ denote the arrow, the only 1\nbd dimensional atom, whose underlying poset is $I = \set{0^- < 1 > 0^+}$.
Given a graded poset $P$ and a closed subset $K \subseteq P$, the \emph{partial cylinder on $P$ relative to $K$} is the graded poset $I \pcyl{K} P$ obtained as the pushout
\[\begin{tikzcd}
	{I \times K} & K \\
	{I \times P} & {(I \times P) \coprod_{I \times K} K} 
	\arrow[two heads, from=1-1, to=1-2]
	\arrow[hook', from=1-1, to=2-1]
	\arrow["{(-)}", hook', from=1-2, to=2-2]
	\arrow["q", two heads, from=2-1, to=2-2]
	\arrow["\lrcorner"{anchor=center, pos=0.125, rotate=180}, draw=none, from=2-2, to=1-1]
\end{tikzcd}\]
in $\Pos$.
This is equipped with a canonical projection map $\tau_K\colon I \pcyl{K} P \surj P$.
\end{dfn}

\noindent
Explicitly, an element of $I \pcyl{K} P$ is either
\begin{itemize}
	\item $(x)$ where $x \in K$, or
	\item $(i, x)$ where $i \in I$ and $x \in P \setminus K$,
\end{itemize}
and the partial order is defined by
\begin{align*}
	\faces{}{}(x) & \eqdef \set{(y) \mid y \in \faces{}{}x}, \\
	\faces{}{}(i, x) & \eqdef \begin{cases}
		\set{(0^-, x), (0^+, x)} + \set{(1, y) \mid y \in \faces{}{}x \setminus K} &
		\text{if $i = 1$,} \\
		\set{(i, y) \mid y \in \faces{}{}x \setminus K} + 
		\set{(y) \mid y \in \faces{}{}x \cap K} &
		\text{otherwise}.
	\end{cases}
\end{align*}

\begin{lem} \label{lem:cylinder_projection_is_cartesian}
	Let $P$ be a graded poset and $K \subseteq P$.
	Then the projection $\tau_K\colon I \pcyl{K} P \surj P$ is a cartesian map of posets.
\end{lem}
\begin{proof}
	Let $x \in K$, so $(x) \in I \pcyl{K} P P$.
	Then \( \restr{\tau_K}{\clset{(x)}} \) is an isomorphism, hence cartesian. 
	Else, let \( i \in I \) and \( x \in P \setminus K \), so that $(i, x) \in I \pcyl{K} P P$, and consider \( y \in \faces{}{}x = \faces{}{}\tau_K(i, x) \). 
	If \( y \in K \), then \( (y) \in \faces{}{} (i, x) \) is a lift of \( y \), and since \( \restr{\tau_K}{\clos {(y)}} \) is an isomorphism, it is a cartesian lift.
	Otherwise, if \( y \notin K \), we claim that \( (i, y) \in \faces{}{} (i, x) \) is a cartesian lift of \( y \).
	Indeed, an element in $\clset{(i, x)}$ is either of the form $(z)$ for $z \in K$ or of the form $(j, z)$ for $j \le i$ and $z \in P \setminus K$.
	Then $z \le y$ implies $(z) \le (i, y)$ and $(j, z) \le (i, y)$, respectively.
\end{proof}

\begin{dfn}[Partial Gray cylinder]
	Let $U$ be a regular directed complex and $K \subseteq U$ a closed subset.
	The \emph{partial Gray cylinder on $U$ relative to $K$} is the oriented graded poset $\arr \pcyl{K} U$ whose
\begin{itemize}
	\item underlying graded poset is $I \pcyl{K} P$, and
	\item orientation is specified, for all $\a \in \set{+, -}$, by
\begin{align*}
	\faces{}{\a}(x) & \eqdef \set{(y) \mid y \in \faces{}{\a}x}, \\
	\faces{}{\a}(i, x) & \eqdef \begin{cases}
		\set{(0^\a, x)} + \set{(1, y) \mid y \in \faces{}{-\a}x \setminus K} &
		\text{if $i = 1$,} \\
		\set{(i, y) \mid y \in \faces{}{\a}x \setminus K} + 
		\set{(y) \mid y \in \faces{}{\a}x \cap K} &
		\text{otherwise}.
	\end{cases}
\end{align*}
\end{itemize}
\end{dfn}

\begin{rmk}
	When $K = \varnothing$, the partial Gray cylinder $\arr \pcyl{K} U$ is the Gray product $\arr \gray U$.
	When \( K = U \), it is isomorphic to \( U \).
\end{rmk}

\begin{lem} \label{lem:partial_gray_is_molecule}
	Let $U$ be a molecule and $K \subseteq U$ a closed subset.
	Then 
	\begin{enumerate}
		\item $\arr \pcyl{K} U$ is a molecule, and
		\item $\tau_K\colon \arr \pcyl{K} U \surj U$ is a cartesian map of regular directed complexes.
	\end{enumerate}
	Moreover, if $U$ is round and $K \subseteq \bd{}{}U$, then $\arr \pcyl{K} U$ is round.
\end{lem}
\begin{proof}
	Let $q\colon I \times U \surj I \pcyl{K} P$ be the quotient map appearing in the definition of $I \pcyl{K} P$, and equip $I \times U$ and $I \pcyl{K} P$ with the orientations of $\arr \gray U$ and $\arr \pcyl{K} U$, respectively.
	By \cite[Proposition 7.2.16]{hadzihasanovic2024combinatorics}, $\arr \gray U$ is a molecule.
	We will prove, by induction on submolecules, that for all $J \submol \arr$ and all $V \submol U$, $q(J \gray U)$ is a molecule and a submolecule of $q(\arr \gray U) = \arr \pcyl{K} U$; and that $q(\arr \gray V)$ is round if $V$ is round and $K \cap V \subseteq \bd{}{}V$.
	For all $\a \in \set{+, -}$, $q(\set{0^\a} \gray U)$ is isomorphic to $U$.
	For all $x \in \gr{0}{U}$, $q(\arr \gray \set{x})$ is either a point if $x \in K$, or an arrow if $x \notin K$.
	Next, by inspection, for all $n \in \mathbb{N}$ and $\a \in \set{+, -}$, we have $q(\bd{n}{\a}(\arr \gray U)) = \bd{n}{\a}(\arr \pcyl{K} U)$, which implies that $\bd{n}{\a}(\arr \pcyl{K} U)$ is globular, and
	\[
		\bd{n}{\a}(\arr \gray U) = (\set{0^\a} \gray \bd{n}{\a}U) \cup (\arr \gray \bd{n-1}{-\a}U).
	\]
	Suppose that $U$ is round and $K \subseteq \bd{}{}U$, so in particular $\arr \pcyl{K} U$ has dimension $\dim U + 1$.
	For all $n \le \dim U$,
	\begin{align*}
		\bd{n}{+}&(\arr \pcyl{K} U) \cap \bd{n}{-}(\arr \pcyl{K} U) \\
			& = (q(\set{0^+} \gray \bd{n}{+}U) \cap q(\set{0^-} \gray \bd{n}{-}U)) \cup 
			(q(\set{0^+} \gray \bd{n}{+}U) \cap q(\arr \gray \bd{n-1}{+}U)) \\
			& \qquad \cup (q(\arr \gray \bd{n-1}{-}U) \cap q(\set{0^-} \gray \bd{n}{-}U)) \cup
			(q(\arr \gray \bd{n-1}{-}U) \cap q(\arr \gray \bd{n-1}{+}U)) \\
			& = q(\arr \gray (K \cap \bd{n}{+}U \cap \bd{n}{-}U)) 
			\cup q(\set{0^+} \gray \bd{n-1}{+}U) \\
			& \qquad \cup q(\set{0^-} \gray \bd{n-1}{-}U)
			\cup q(\arr \gray \bd{n-2}{}U) \\
			& = q(\arr \gray (K \cap \bd{n}{+}U \cap \bd{n}{-}U)) \cup \bd{n-1}{}(\arr \pcyl{K} U).
	\end{align*}
	The first component is included in the second one because $K \subseteq \bd{}{}U$ when $n = \dim U$, and because $U$ is round when $n < \dim U$, so we conclude that $\arr \pcyl{K} U$ is round under the assumptions.
	If $U$ is an atom, then either $K = U$ and $\arr \pcyl{K} U$ is isomorphic to $U$, or $K \subseteq \bd{}{}U$ and $U$ is round, so the previous argument together with the inductive hypothesis suffice to conclude that $\arr \pcyl{K} U$ is an atom.
	If $U$ is not an atom, it splits into $V \cp{k} W$, and by \cite[Proposition 7.2.16]{hadzihasanovic2024combinatorics} $\arr \gray U$ splits into $(\arr \gray W) \gencp{k+1} (\arr \gray V)$.
	By inspection, $q$ preserves this generalised pasting, and by the inductive hypothesis we conclude that $\arr \pcyl{K} U$ is a molecule.
	Finally, $\tau_K$ is cartesian by Lemma \ref{lem:cylinder_projection_is_cartesian}, and the properties of $q$, along with the fact that the projection $\arr \gray U \surj U$ is a map of regular directed complexes, allow us to conclude that it is a cartesian map of regular directed complexes.
\end{proof}

\begin{dfn}[Unit]
	Let $u\colon U \to X$ be a pasting diagram in a diagrammatic set.
	The \emph{unit on $u$} is the degenerate pasting diagram $\un u\colon u \celto u$ defined by $u \after \tau_{\bd{}{}U}\colon \arr \pcyl{\bd{}{}U} U \to X$.
\end{dfn}

\begin{dfn}[Left unitor]
	Let $u\colon U \to X$ be a pasting diagram in a diagrammatic set and let $\iota\colon v \submol \bd{}{-}u$ be a rewritable subdiagram of shape $V$ in its input boundary.
	Let $K \eqdef \bd{}{}U \setminus \inter \iota(V)$.
	The \emph{left unitor of $u$ at $\iota$} is the degenerate pasting diagram $\lun{\iota}u\colon u \celto \un v \cpsub{\iota} u$ defined by $u \after \tau_{K}\colon \arr \pcyl{K} U \to X$.
\end{dfn}

\begin{dfn}[Right unitor]
	Let $u\colon U \to X$ be a pasting diagram in a diagrammatic set and let $\iota\colon v \submol \bd{}{+}u$ be a rewritable subdiagram of shape $V$ in its output boundary.
	Let $K \eqdef \bd{}{}U \setminus \inter \iota(V)$.
	The \emph{right unitor of $u$ at $\iota$} is the degenerate pasting diagram $\run{\iota}u\colon u \subcp{\iota} \un v \celto u$ defined by $u \after \tau_{K}\colon \arr \pcyl{K} U \to X$.
\end{dfn}

\noindent We will simply write $\lun{} u$ and $\run{} u$ when $\iota$ is an isomorphism.

\begin{rmk}
	If $u$ is round, then by Lemma \ref{lem:partial_gray_is_molecule} so are $\un u$, $\lun{\iota}u$, and $\run{\iota}u$.
\end{rmk}

\begin{dfn}[Inverted partial Gray cylinder]
	Let $U$ be a molecule, $n \eqdef \dim{U}$, and $K \subseteq \bd{}{+}U$ a closed subset.
	The \emph{left-inverted partial Gray cylinder on $U$ relative to $K$} is the oriented graded poset $\lcyl{K} U$ whose
\begin{itemize}
	\item underlying graded poset is $I \pcyl{K} P$, and
	\item orientation is as in $\arr \pcyl{K} U$, except for all $x \in \gr{n}{U}$ and $\a \in \set{+, -}$
\begin{align*}
	\faces{}{-}(1, x) &\eqdef \set{(0^-, x), (0^+, x)} + \set{(1, y) \mid y \in \faces{}{+}x \setminus K}, \\
	\faces{}{+}(1, x) &\eqdef \set{(1, y) \mid y \in \faces{}{-}x}, \\
	\faces{}{\a}(0^+, x) &\eqdef \set{(0^+, y) \mid y \in \faces{}{-\a}x \setminus K} + 
		\set{(y) \mid y \in \faces{}{-\a}x \cap K}.
\end{align*}
\end{itemize}
	Dually, if $K \subseteq \bd{}{-}U$, the \emph{right-inverted partial Gray cylinder on $U$ relative to $K$} is the oriented graded poset $\rcyl{K}{U}$ whose
\begin{itemize}
	\item underlying graded poset is $I \pcyl{K} P$, and
	\item orientation is as in $\arr \pcyl{K} U$, except for all $x \in \gr{n}{U}$ and $\a \in \set{+, -}$
\begin{align*}
	\faces{}{-}(1, x) &\eqdef \set{(1, y) \mid y \in \faces{}{+}x}, \\
	\faces{}{+}(1, x) &\eqdef \set{(0^-, x), (0^+, x)} + \set{(1, y) \mid y \in \faces{}{-}x \setminus K}, \\
	\faces{}{\a}(0^-, x) &\eqdef \set{(0^-, y) \mid y \in \faces{}{-\a}x \setminus K} + 
		\set{(y) \mid y \in \faces{}{-\a}x \cap K}.
\end{align*}
\end{itemize}
\end{dfn}

\begin{lem} \label{lem:inverted_cylinders_are_molecules}
	Let $U$ be a molecule, let $K \subseteq \bd{}{+}U$ and $K' \subseteq \bd{}{-}U$ be closed subsets, and let $p\colon U \to V$ be a cartesian map of regular directed complexes such that $\dim{V} < \dim{U}$.
	Then
	\begin{enumerate}
		\item $\lcyl{K}{U}$ and $\rcyl{K'}{U}$ are molecules, and
		\item $p \after \tau_K\colon \lcyl{K}{U} \to V$ and $p \after \tau_{K'}\colon \rcyl{K'}{U} \to V$ are cartesian maps of regular directed complexes.
	\end{enumerate}
	Moreover, if $U$ is round, then $\lcyl{K}{U}$ and $\rcyl{K'}{U}$ are round.
\end{lem}
\begin{proof}
	We will prove the statement for $\lcyl{K}{U}$, the case of $\rcyl{K'}{U}$ being dual.
	By construction, $\bd{}{-}\lcyl{K}{U}$ and $\bd{}{+}\lcyl{K}{U}$ are, respectively, of the form
	\[
		U \cpsub{} (\arr \pcyl{K} \bd{}{+}U) \subcp{} \dual{n}{U}
		\quad \text{ and } \quad 
		\arr \pcyl{K \cap \bd{}{-}U} \bd{}{-}U.
	\]
	By Lemma \ref{lem:partial_gray_is_molecule}, both are molecules.
	Moreover, if $U$ is round, since $K \subseteq \bd{}{+}U$, we have $K \cap \bd{}{-}U \subseteq \bd{}{}(\bd{}{-}U)$, so $\bd{}{+}\lcyl{K}{U}$ is round.
	Consequently, $\bd{}{-}\lcyl{K}{U}$ is also round.
	We now proceed by induction on the layering dimension of $U$.
	Let $n \eqdef \dim{U}$, which is necessarily $> 0$ by the fact that $p$ exists.
	If $\lydim{U} = -1$, then $U$ is an atom, and the previous argument suffices to prove that $\lcyl{K}{U}$ is an atom.
	If $\lydim{U} < n-1$, then $\gr{n}{U} = \set{x}$, and $K_x \eqdef K \cap \clset{x} \subseteq \bd{}{+}x$ by \cite[Lemma 4.3.14]{hadzihasanovic2024combinatorics}.
	Then $\lcyl{K_x}{\clset{x}}$ is well-defined as an atom, and $\bd{}{-}\lcyl{K_x}{\clset{x}} \submol \bd{}{-}\lcyl{K}{U}$.
	It follows that $\lcyl{K}{U} = \bd{}{-}\lcyl{K}{U} \subcp{n-1} \lcyl{K_x}{\clset{x}}$ is well-defined as a molecule.
	Finally, suppose that $\lydim{U} = n-1$, and pick an $(n-1)$\nbd layering $(\order{i}{U})_{i=1}^m$ of $U$.
	By \cite[Lemma 4.1.6]{hadzihasanovic2024combinatorics}, we have $K_i \eqdef K \cap \order{i}{U} \subseteq \bd{}{+}\order{i}{U}$ for all $i \in \set{1, \ldots, m}$, so $\lcyl{K_i}{\order{i}{U}}$ is well-defined as a molecule, and so is
	\[
		\lcyl{K}{U} = ((\bd{}{-}\lcyl{K}{U} \subcp{n-1} \lcyl{K_m}{\order{m}{U}}) \subcp{} \lcyl{K_{m-1}}{\order{m-1}{U}} \ldots ) \subcp{} \lcyl{K_1}{\order{1}{U}}.
	\]
	It is straightforward to show that this is round when $U$ is round.
	Finally, $p \after \tau_K$ is cartesian by Lemma \ref{lem:cylinder_projection_is_cartesian}, so it only remains to show that it is a map of regular directed complexes.
	This follows from Lemma \ref{lem:partial_gray_is_molecule} on the closure of every element which is not of the form $(0^+, x)$ or $(1, x)$ for some $x \in \gr{n}{U}$.
	Let $x \in \gr{n}{U}$.
	Then $\restr{(p \after \tau_K)}{\clset{(0^+, x)}}$ is equal up to isomorphism to $\rev{(\restr{p}{\clset{x}})}$, which we already know to be a map of regular directed complexes.
	It only remains to show that, for all $k \leq n$ and $\a \in \set{+, -}$, $p(\tau_K(\bd{k}{\a}(1, x))) = \bd{k}{\a}p(x)$, and that $\restr{(p \after \tau_K)}{\bd{k}{\a}(1,x)}$ is final onto its image.
	We have
	\[
		\tau_K(\bd{k}{\a}(1, x)) =
		\begin{cases}
			\clset{x} 
				& \text{if $k = n$ and $\a = -$}, \\
			\bd{}{-}x
				& \text{if $k = n$ and $\a = +$, or if $k = n-1$}, \\
			\bd{k}{\a}x
				& \text{otherwise}.
		\end{cases}
	\]
	Since $\dim p(x) < n$, we have $p(\bd{k}{\a}x) = \bd{k}{\a}p(x) = \clset{p(x)}$ for all $k \geq n-1$ and $\a \in \set{+, -}$, which proves that $p \after \tau_K$ is compatible with boundaries.
	Moreover, finality of $\restr{p}{\bd{k}{\a}x}$ onto its image, together with the fact that zig-zags in a closed $W \subseteq U$ can be lifted to zig-zags in $\arr \pcyl{K \cap W} U$, imply finality of $\restr{(p \after \tau_K)}{\bd{k}{\a}(1, x)}$ in all cases except when $k = n$ and $\a = -$.
	In this last case, finality of $\restr{p}{\bd{}{+}x}$ onto its image takes care of all identified pairs of elements except the pair of $(0^-, x)$ and $(0^+, x)$.
	By \cite[Lemma 6.2.4]{hadzihasanovic2024combinatorics}, there exists $y \in \faces{}{+}x$ such that $p(y) = p(x)$, which implies the existence of a zig-zag between $(0^-, x)$ and $(0^+, x)$ in $\bd{}{-}(1, x)$ all mapped to $p(x)$ by $p \after \tau_K$.
	This concludes the proof.
\end{proof}

\begin{dfn}[Left invertor]
	Let $u\colon U \to X$ be a degenerate pasting diagram in a diagrammatic set.
	The \emph{left invertor of $u$} is the degenerate pasting diagram $\linv u\colon u \cp{} \rev{u} \celto \un {(\bd{}{-}u)}$ defined by $u \after \tau_{\bd{}{+}U}\colon \lcyl{\bd{}{+}U} U \to X$.
\end{dfn}

\begin{dfn}[Right invertor]
	Let $u\colon U \to X$ be a degenerate pasting diagram in a diagrammatic set.
	The \emph{right invertor of $u$} is the degenerate pasting diagram $\rinv u\colon \un {(\bd{}{+}u)} \celto \rev{u} \cp{} u$ defined by $u \after \tau_{\bd{}{-}U}\colon \rcyl{\bd{}{-}U} U \to X$.
\end{dfn}

\begin{rmk}
	By Lemma \ref{lem:inverted_cylinders_are_molecules}, if $u$ is a degenerate round diagram, then $\linv u$ and $\rinv u$ are also round.
\end{rmk}

\noindent We conclude this section with a couple of pictures of the families of degenerate diagrams we just introduced. 
Suppose \( f \colon u \celto v \) is a \( 1 \)\nbd cell in a diagrammatic set \( X \).
Then, the degenerate \( 2 \)\nbd cells \( \lun{}f \colon \arr \pcyl{\bd{}{+} \arr} \arr \to X \) and \( \run{}f \colon \arr \pcyl{\bd{}{-} \arr} \arr \to X \) are respectively depicted as
\begin{center}
	\begin{tikzcd}
		& u &&& u && {v.} \\
		u && v & {\text{and}} && v
		\arrow["f", curve={height=-6pt}, from=1-2, to=2-3]
		\arrow[""{name=0, anchor=center, inner sep=0}, "f", curve={height=-12pt}, from=1-5, to=1-7]
		\arrow["f"', curve={height=6pt}, from=1-5, to=2-6]
		\arrow["{\un u}", curve={height=-6pt}, from=2-1, to=1-2]
		\arrow[""{name=1, anchor=center, inner sep=0}, "f"', curve={height=12pt}, from=2-1, to=2-3]
		\arrow["{\un v}"', curve={height=6pt}, from=2-6, to=1-7]
		\arrow["{\lun f}"{description, pos=0.6}, between={0.2}{1}, Rightarrow, from=1, to=1-2]
		\arrow["{\run f}"{description, pos=0.4}, between={0}{0.8}, Rightarrow, from=2-6, to=0]
	\end{tikzcd}
\end{center}
Since \( \lun{}f \) is degenerate, we have the left invertor \( \linv (\lun{} f) \colon (\lun{} f) \cp{} \rev{(\lun{} f)} \celto \un f \).
This is a \( 3 \)\nbd cell whose input and output \( 2 \)\nbd boundaries are respectively
\begin{center}
	\begin{tikzcd}[column sep=scriptsize,row sep=large]
		u && u && v & {\text{and}} & u &&&& v.
		\arrow["{\un u}", from=1-1, to=1-3]
		\arrow[""{name=0, anchor=center, inner sep=0}, "f"', curve={height=30pt}, from=1-1, to=1-5]
		\arrow[""{name=1, anchor=center, inner sep=0}, "f", curve={height=-30pt}, from=1-1, to=1-5]
		\arrow["f", from=1-3, to=1-5]
		\arrow[""{name=2, anchor=center, inner sep=0}, "f"', curve={height=30pt}, from=1-7, to=1-11]
		\arrow[""{name=3, anchor=center, inner sep=0}, "f", curve={height=-30pt}, from=1-7, to=1-11]
		\arrow["{\lun f}"{pos=0.6}, between={0.2}{1}, Rightarrow, from=0, to=1-3]
		\arrow["{\rev{(\lun f)}}", between={0}{0.8}, Rightarrow, from=1-3, to=1]
		\arrow["{\un f}", between={0.2}{0.8}, Rightarrow, from=2, to=3]
	\end{tikzcd}
\end{center}

%% file: equivalences.tex
\section{Equivalences} \label{sec:equivalences}

\subsection{Definition and closure properties}

\begin{dfn}[Lax and colax solutions to equations] \label{dfn:lax_colax_solution}
	Let $X$ be a diagrammatic set and let $\Phi\colon A \to \rd{X}$ be a parametrised family of round diagrams.
	Each $v \in \rd{X}$ determines an \emph{equation $\Phi(x) \qeq v$ in the indeterminate $x \in A$}.
	A \emph{lax solution} to $\Phi(x) \qeq v$ is a pair of
	\begin{enumerate}
		\item $a \in A$ such that $\Phi(a)$ is parallel to $v$, and
		\item a round diagram $h\colon \Phi(a) \celto v$.
	\end{enumerate}
	Dually, a \emph{colax solution} is pair of $a \in A$ and a round diagram $h\colon v \celto \Phi(a)$.
\end{dfn}

\begin{dfn}[Equivalence in a diagrammatic set]
	Let $e$ be a round diagram in a diagrammatic set $X$, $n \eqdef \dim{e} > 0$.
	We say that $e$ is an \emph{equivalence} if, for all parallel $v, w \in \gr{n-1}{\rd X}$,
	\begin{enumerate}
		\item for all rewritable subdiagrams $\iota\colon \bd{}{+}e \submol v$ and for all round diagrams \( u \colon \subs{v}{\bd{}{-}e}{\bd{}{+}e} \celto w \), there exist a round diagram \( u' \colon v \celto w \) and a round diagram $h\colon e \cpsub{\iota} u' \celto u$ which is an equivalence, and
		\item for all rewritable subdiagrams $\iota\colon \bd{}{-}e \submol w$ and for all round diagrams \( u \colon v \celto \subs{w}{\bd{}{+}e}{\bd{}{-}e} \), there exist a round diagram \( u' \colon v \celto w \) and a round diagram $h\colon u' \subcp{\iota} e \celto u$ which is an equivalence.
	\end{enumerate}
	We write $\eqv X$ for the set of equivalences in $X$.
\end{dfn}

\begin{rmk} \label{rmk:equations_and_subdiagrams}
	Using the terminology of Definition \ref{dfn:lax_colax_solution}, $e$ is an \emph{equivalence} if, for all parallel $v, w \in \gr{n-1}{\rd X}$,
	\begin{enumerate}
		\item for all rewritable subdiagrams $\iota\colon \bd{}{+}e \submol v$, every well-formed equation $e \cpsub{\iota} x \qeq u$ in the indeterminate $x \in \gr{n}{\rd X(v, w)}$ has a lax solution $h\colon e \cpsub{\iota} u' \celto u$ such that $h$ is an equivalence, and
		\item for all rewritable subdiagrams $\iota\colon \bd{}{-}e \submol w$, every well-formed equation $x \subcp{\iota} e \qeq u$ in the indeterminate $x \in \gr{n}{\rd X(v, w)}$ has a lax solution $h\colon u' \subcp{\iota} e \celto u$ such that $h$ is an equivalence.
	\end{enumerate}
	By \emph{well-formed equation}, we mean that replacing $x$ with any $u' \in \gr{n}{\rd X(v, w)}$ results in a round diagram parallel to $u$.
	Given $e, u \in \gr{n}{\rd X}$, well-formed equations $e \cpsub{} x \qeq u$ are in bijection with subdiagrams $\bd{}{-}e \submol \bd{}{-}u$.
	Indeed, given $\iota\colon \bd{}{+}e \submol v$ and a round diagram $u'\colon v \celto w$, the round diagram $e \cpsub{\iota} u'$ has type $\subs{v}{\bd{}{-}e}{\iota(\bd{}{+}e)} \celto w$, which contains $\bd{}{-}e$ as a subdiagram of its input.
	Conversely, given $j\colon \bd{}{-}e \submol \bd{}{-}u$, we let $v \eqdef \subs{\bd{}{-}u}{\bd{}{+}e}{j(\bd{}{-}e)}$ and $\iota\colon \bd{}{+}e \submol v$ be the evident subdiagram.
	Dually, well-formed equations $x \subcp{} e \qeq u$ are in bijection with subdiagrams $\bd{}{+}e \submol \bd{}{+}u$.
\end{rmk}

\begin{comm}
	The definition of equivalence is coinductive; we make it more explicit for those unfamiliar with this style.
	Given a set $A \subseteq \rd X$, we let $\E(A)$ be the set of round diagrams $e$ such that, letting $n \eqdef \dim{e} > 0$, for all parallel $v, w \in \gr{n-1}{\rd X}$,
	\begin{enumerate}
		\item for all rewritable subdiagrams $\iota\colon \bd{}{+}e \submol v$, every well-formed equation $e \cpsub{\iota} x \qeq u$ in the indeterminate $x \in \gr{n}{\rd X(v, w)}$ has a lax solution $h\colon e \cpsub{\iota} u' \celto u$ such that $h \in A$, and
		\item for all rewritable subdiagrams $\iota\colon \bd{}{-}e \submol w$, every well-formed equation $x \subcp{\iota} e \qeq u$ in the indeterminate $x \in \gr{n}{\rd X(v, w)}$ has a lax solution $h\colon u' \subcp{\iota} e \celto u$ such that $h \in A$.
	\end{enumerate}
	Then $\E$ defines an order-preserving operator on the power set $\powerset{(\rd X)}$, which by the Knaster--Tarski theorem admits a greatest fixed point
	\[
		\eqv X = \bigcap_{k \geq 0} \E^k(\rd X),
	\]
	and this is the set of equivalences in $X$.
	This definition comes with the following proof method: given any $A \subseteq \rd X$, if $A \subseteq \E(A)$, then $A \subseteq \eqv X$.
\end{comm}

\begin{rmk}
	The definition of equivalence may seem biased towards lax, rather than colax solutions, but this is illusory: we will find that requiring colax solutions results in the same notion.
\end{rmk} 

\begin{exm}
	A \( 1 \)\nbd dimensional round diagram \( e \colon x \celto y \) is is an equivalence if, for all \( 0 \)\nbd cells \( z \) and all round diagrams \( u \colon x \celto z \), there exists a round diagram \( u' \colon y \celto z \), together with an equivalence \( h \colon e \cp{} u' \celto u \), and dually, for all round diagrams \( u \colon z \celto y \), there exists a round diagram \( u' \colon z \celto x \), together with an equivalence \( h \colon u' \cp{} e \celto u \).
	The informal picture to keep in mind is the following:
	\begin{center}
		\begin{tikzcd}
			x && z & {\text{and dually}} & z && y. \\
			& y &&&& x
			\arrow[""{name=0, anchor=center, inner sep=0}, "{\forall u}"{description}, curve={height=-12pt}, from=1-1, to=1-3]
			\arrow["e"{description}, curve={height=6pt}, from=1-1, to=2-2]
			\arrow[""{name=1, anchor=center, inner sep=0}, "{\forall u}"{description}, curve={height=-12pt}, from=1-5, to=1-7]
			\arrow["{\exists u'}"{description}, curve={height=6pt}, from=1-5, to=2-6]
			\arrow["{\exists u'}"{description}, curve={height=6pt}, from=2-2, to=1-3]
			\arrow["e"{description}, curve={height=6pt}, from=2-6, to=1-7]
			\arrow["{\exists h}"{description}, between={0}{0.8}, Rightarrow, from=2-2, to=0]
			\arrow["{\exists h}"{description}, between={0}{0.8}, Rightarrow, from=2-6, to=1]
		\end{tikzcd}
	\end{center} 
\end{exm}

\noindent Given $A \subseteq \rd X$, we let $\T(A)$ denote the closure of $A$ under the following clauses: for all $n \in \mathbb{N}$ and all $h \in \gr{n+1}{\rd{X}}$ of type $u \cpsub{} v \celto w$ or $v \subcp{} u \celto w$ with $u, v, w \in \gr{n}{\rd{X}}$,
\begin{enumerate}
	\item if $h, u, v \in \T(A)$, then $w \in \T(A)$, and
	\item if $h, u, w \in \T(A)$, then $v \in \T(A)$.
\end{enumerate}
Intuitively, $\T(A)$ is the closure of $A$ under ``composition'' and ``division'' of $n$\nbd dimensional round diagrams as witnessed by $(n+1)$\nbd dimensional round diagrams; this is a form of 2-out-of-3 property.

Let $\dgn X \eqdef \set{ u \in \rd X \mid \text{$u$ is degenerate} }$.
Our next goal is to prove that $\eqv X = \T(\eqv X \cup \dgn X)$, that is, equivalences include all degenerate round diagrams and are closed under 2-out-of-3.

\begin{lem} \label{lem:weak_inverse_in_TA}
	Let $X$ be a diagrammatic set, $A \subseteq \rd X$, and let $e \in A \cap \E(A)$ be of type $u \celto v$.
	Then there exists $e^* \in \T(A)$ of type $v \celto u$.
\end{lem}
\begin{proof}
	The equation $e \cp{} x \qeq e$ is well-defined and, since $e \in \E(A)$, it admits a lax solution $h\colon e \cp{} u \celto e$ with $h \in A$.
	Since $h, e \in A$, it follows that $u \in \T(A)$.
	Now, $u$ has type $v \celto v$, so the equation $x \cp{} e \qeq u$ is well-defined and admits a lax solution $k\colon e^* \cp{} e \celto u$ with $k \in A$.
	Then $e^* \in \T(A)$ has type $v \celto u$.
\end{proof}

\begin{lem} \label{lem:weak_inverse_in_Teqv}
	Let $X$ be a diagrammatic set and let $e \in \eqv X$ have type $u \celto v$.
	Then there exists $e^* \in \T(\eqv X)$ of type $v \celto u$.
\end{lem}
\begin{proof}
	Follows from Lemma \ref{lem:weak_inverse_in_TA} with $A \eqdef \eqv X$.
\end{proof}

\begin{lem} \label{lem:composition_in_TA}
	Let $X$ be a diagrammatic set, $\dgn X \subseteq A \subseteq \rd X$, $n \in \mathbb{N}$, and $u, v \in \gr{n}{A}$.
	Then
	\begin{enumerate}
		\item if a pasting $u \cpsub{} v$ is defined, then $u \cpsub{} v \in \T(A)$,
		\item if a pasting $v \subcp{} u$ is defined, then $v \subcp{} u \in \T(A)$.
	\end{enumerate}
\end{lem}
\begin{proof}
	Suppose $u \cpsub{} v$ is defined and consider the unit $\un (u \cpsub{} v)\colon u \cpsub{} v \celto u \cpsub{} v$.
	Then $u, v \in A$ and $\un(u \cpsub{} v) \in \dgn X \subseteq A$.
	It follows that $u \cpsub{} v \in \T(A)$.
	The case where $v \subcp{} u$ is defined is dual.
\end{proof}

\begin{dfn}[Unbiased set of solutions]
	Let $X$ be a diagrammatic set, $A \subseteq \rd X$ and $e \in \E(A)$.
	We say that $A$ is \emph{unbiased for $e$} if equations $e \cpsub{} x \qeq u$ admit a pair of a lax solution $h\colon e \cpsub{} u' \celto u$ and a colax solution $h^*\colon u \celto e \cpsub{} u'$ with $h, h^* \in A$, and similarly for equations $x \subcp{} e \qeq u$.
	Given $B \subseteq \E(A)$, we say that $A$ is \emph{unbiased for $B$} if it is unbiased for all $e \in B$.
\end{dfn}

\begin{lem} \label{lem:Teqv_unbiased_for_eqv}
	Let $X$ be a diagrammatic set. 
	Then $\eqv X \subseteq \E(\T(\eqv X))$, and $\T(\eqv X)$ is unbiased for $\eqv X$.
\end{lem}
\begin{proof}
	First of all, $\eqv X = \E(\eqv X) \subseteq \E(\T(\eqv X))$ because $\E$ is order-preserving.
	Let $e \in \eqv X$, and consider an equation $e \cpsub{} x \qeq u$.
	By definition, this admits a lax solution $h\colon e \cpsub{} u' \celto u$ with $h \in \eqv X$.
	By Lemma \ref{lem:weak_inverse_in_Teqv}, there also exists a colax solution $h^*\colon u \celto e \cpsub{} u'$ in $\T(\eqv X)$.
	The case of equations $x \subcp{} e \qeq u$ is dual.
\end{proof}

\begin{lem} \label{lem:Tdgn_unbiased_for_dgn}
	Let $X$ be a diagrammatic set.
	Then $\dgn X \subseteq \E(\T(\dgn X))$, and $\T(\dgn X)$ is unbiased for $\dgn X$.
\end{lem}
\begin{proof}
	Let $e \in \dgn X$ have type $v \celto w$, and consider an equation $e \cpsub{\iota} x \qeq u$, which by Remark \ref{rmk:equations_and_subdiagrams} corresponds to a subdiagram $j\colon v \submol \bd{}{-}u$.
	The left invertor $\linv{e}\colon e \cp{} \rev{e} \celto \un v$ and reverse left unitor $\rev{(\lun{j} u)}\colon \un v \cpsub{j} u \celto u$ are both in $\dgn X$.
	Moreover, there is an evident pasting $h \eqdef \linv{e} \cpsub{} \rev{(\lun{j} u)}$, which is a round diagram of type $e \cpsub{\iota} (\rev{e} \cpsub{j} u) \celto u$, that is, a lax solution for $e \cpsub{\iota} x \qeq u$.
	By Lemma \ref{lem:composition_in_TA}, $h \in \T(\dgn X)$.
	Dually, $h^* \eqdef \lun{j} u \subcp{} \rev{(\linv{e})}$ is a colax solution of type $u \celto e \cpsub{\iota} (\rev{e} \cpsub{j} u)$.
	The case of equations $x \subcp{\iota} e \qeq u$ is dual, using right invertors and right unitors.
\end{proof}

\begin{lem} \label{lem:unbiased_main_lemma}
	Let $X$ be a diagrammatic set and $\dgn X \subseteq A \subseteq \rd X$.
	Suppose that $A \subseteq \E(\T(A))$ and $\T(A)$ is unbiased for $A$.
	Then $\T(A) \subseteq \eqv X$.
\end{lem}
\begin{proof}
	We will prove, by structural induction on the definition of \( \T(A) \), that for all $e \in \T(A)$, we have $e \in \E(\T(A))$ and $\T(A)$ is unbiased for $e$.
	This will prove that $\T(A) \subseteq \E(\T(A))$, so by coinduction $\T(A) \subseteq \eqv X$.
	It will suffice to consider equations of the form $e \cpsub{} x \qeq u$, the other kind being dual.

	The base case $e \in A$ holds by assumption.
	Next, consider a round diagram $k\colon v \cpsub{} w \celto e$, and assume the inductive hypothesis of $k, v, w \in \T(A)$.
	By Remark \ref{rmk:equations_and_subdiagrams}, since $\bd{}{-}v \submol \bd{}{-}e \submol \bd{}{-}u$, there is a well-formed equation $v \cpsub{} x \qeq u$, which by the inductive hypothesis admits both a lax and a colax solution
	\[ k_v\colon v \cpsub{} u_v \celto u, \quad \quad k_v^*\colon u \celto v \cpsub{} u_v \]
	in $\T(A)$.
	Next, since $\bd{}{-}w = \subs{\bd{}{-}e}{\bd{}{+}v}{\bd{}{-}v} \submol \subs{\bd{}{-}u}{\bd{}{+}v}{\bd{}{-}v} = \bd{}{-}u_v$, we have a well-formed equation $w \cpsub{} x \qeq u_v$, which by the inductive hypothesis admits a pair of a lax and a colax solution
	\[ k_w\colon w \cpsub{} u_w \celto u_v, \quad \quad k_w^*\colon u_v \celto w \cpsub{} u_w \]
	in $\T(A)$.
	Now, the evident pastings 
	\[
		k_w \cpsub{} k_v\colon (v \cpsub{} w) \cpsub{} u_w \celto u,
		\quad \quad
		k_v^* \subcp{} k_w^*\colon u \celto (v \cpsub{} w) \cpsub{} u_w
	\]
	are defined and belong to $\T(A)$ by Lemma \ref{lem:composition_in_TA}.
	Since $k \in \T(A) \cap \E(\T(A))$, by Lemma \ref{lem:weak_inverse_in_TA} there exists $k^*\colon e \celto v \cpsub{} w$ in $\T(A)$.
	Then the evident pastings
	\[
		k^* \cpsub{} (k_w \cpsub{} k_v)\colon e \cpsub{} u_w \celto u,
		\quad \quad
		(k_v^* \subcp{} k_w^*) \subcp{} k\colon u \celto e \cpsub{} u_w
	\]
	are defined and belong to $\T(A)$.
	These exhibit a pair of a lax and colax solution for $e \cpsub{} x \qeq u$.
	The case of $k\colon v \subcp{} w \celto e$ is analogous.
	
	Next, consider a round diagram $k\colon v \cpsub{\iota} e \celto w$, assuming the inductive hypothesis of $k, v, w \in \T(A)$.
	Let $z \eqdef \bd{}{+}v$ and consider the equation $x \cp{} v \qeq \un z$.
	This admits a pair of a lax and colax solution 
	\[
		h\colon v^* \cp{} v \celto \un z,
		\quad \quad
		h^*\colon \un z \celto v^* \cp{} v
	\]
	in $\T(A)$.
	Moreover, we have a subdiagram $j\colon z \submol \bd{}{-}u$ obtained by composing $\iota\colon z \submol \bd{}{-}e$ with $\bd{}{-}e \submol \bd{}{-}u$.
	We let
	\[
		k_1 \eqdef h \cpsub{} \rev{(\lun{j} u)} \colon v^* \cpsub{} (v \cpsub{j} u) \celto u,
		\quad \quad
		k_1^* \eqdef \lun{j}u \subcp{} h^*\colon u \celto v^* \cpsub{} (v \cpsub{j} u) 
	\]
	be the evident pastings, which are in $\T(A)$ by Lemma \ref{lem:composition_in_TA}.
	Now, by construction $\bd{}{-} (v \cpsub{j} u) = \subs{\bd{}{-}u}{\bd{}{-}v}{j(z)}$ and $\bd{}{-} w = \subs{\bd{}{-}e}{\bd{}{-}v}{\iota(z)}$, so the subdiagram $\bd{}{-}e \submol \bd{}{-}u$ induces a subdiagram $\bd{}{-}w \submol \bd{}{-} (v \cpsub{j} u)$ and we have a well-formed equation $w \cpsub{} x \qeq v \cpsub{j} u$.
	This admits a pair of a lax and colax solution
	\[
		k_w\colon w \cpsub{} u_w \celto v \cpsub{j} u, 
		\quad \quad
		k_w^*\colon v \cpsub{j} u \celto w \cpsub{} u_w
	\]
	in $\T(A)$, and we let
	\[
		k_2 \eqdef k_w \cpsub{} k_1\colon v^* \cpsub{} (w \cpsub{} u_w) \celto u,
		\quad \quad
		k_2^* \eqdef k_1^* \subcp{} k_w^*\colon u \celto v^* \cpsub{} (w \cpsub{} u_w)
	\]
	be the evident pastings, which are in $\T(A)$ by Lemma \ref{lem:composition_in_TA}.
	Now, since $k \in \T(A) \cap \E(\T(A))$, by Lemma \ref{lem:weak_inverse_in_TA} there exists $k^*\colon w \celto v \cpsub{\iota} e$ in $\T(A)$.
	Then we define the evident pastings
	\[
		k_3 \eqdef k \cpsub{} k_2\colon ((v^* \cp{} v) \cpsub{\iota} e) \cpsub{} u_w \celto u,
		\quad
		k_3^* \eqdef k_2^* \subcp{} k^*\colon u \celto ((v^* \cp{} v) \cpsub{\iota} e) \cpsub{} u_w,
	\]
	followed by the evident pastings
	\[
		k_4 \eqdef h^* \cpsub{} k_3\colon (\un z \cpsub{\iota} e) \cpsub{} u_w \celto u,
		\quad \quad
		k_4^* \eqdef k_3^* \subcp{} h\colon u \celto (\un z \cpsub{\iota} e) \cpsub{} u_w,
	\]
	as well as the evident pastings
	\[
		k_5 \eqdef \lun{\iota} e \cpsub{} k_4\colon e \cpsub{} u_w \celto u,
		\quad \quad
		k_5^* \eqdef k_4^* \subcp{} \rev{(\lun{\iota} e)}\colon u \celto e \cpsub{} u_w,
	\]
	all of which are in $\T(A)$ by repeated applications of Lemma \ref{lem:composition_in_TA}.
	These exhibit a pair of a lax and colax solution to $e \cpsub{} x \qeq u$.

	Finally, consider a round diagram $k\colon e \subcp{} v \celto w$, assuming the inductive hypothesis of $k, v, w \in \T(A)$.
	We have $\bd{}{-}w = \bd{}{-}e \submol \bd{}{-}u$, so the equation $w \cpsub{} x \qeq u$ is well-formed and admits a pair of a lax and colax solution
	\[
		k_w\colon w \cpsub{} u_w \celto u, 
		\quad \quad
		k_w^*\colon u \celto w \cpsub{} u_w
	\]
	in $\T(A)$.
	Moreover, since $k \in \T(A) \cap \E(\T(A))$, by Lemma \ref{lem:weak_inverse_in_TA} there exists $k^*\colon w \celto e \subcp{} v$ in $\T(A)$.
	Then the evident pastings
	\[
		k \cpsub{} k_w\colon e \cpsub{} (v \cpsub{} u_w) \celto u,
		\quad \quad
		k_w^* \subcp{} k^*\colon u \celto e \cpsub{} (v \cpsub{} u_w)
	\]
	are both in $\T(A)$ by Lemma \ref{lem:composition_in_TA}, and exhibit a pair of a lax and colax solution for $e \cpsub{} x \qeq u$.
	This completes the inductive step and the proof.
\end{proof}

\begin{thm} \label{thm:equivalences_2_out_of_3}
	Let $X$ be a diagrammatic set.
	Then
	\begin{enumerate}
		\item $\dgn X \subseteq \eqv X$,
		\item $\eqv X = \T(\eqv X)$.
	\end{enumerate}
\end{thm}
\begin{proof}
	Let $A \eqdef \eqv X \cup \dgn X$.
	Then obviously $\eqv X \subseteq A \subseteq \T(A)$.
	Conversely, by Lemma \ref{lem:Teqv_unbiased_for_eqv} and Lemma \ref{lem:Tdgn_unbiased_for_dgn}, we have $A \subseteq \E(\T(A))$ and $\T(A)$ is unbiased for $A$.
	By Lemma \ref{lem:unbiased_main_lemma}, we conclude that $\T(A) \subseteq \eqv X$.
\end{proof}

\begin{cor} \label{cor:2_out_of_3_consequences}
	Let $X$ be a diagrammatic set, $e \in \eqv X$, and $h \in \rd X$ with $\dim{e} = \dim{h}$.
	Then
	\begin{enumerate}
		\item if $e$ has type $u \celto v$, then there exists $e^* \in \eqv X$ of type $v \celto u$,\label{2_out_of_3_inv}
		\item if a pasting $e \cpsub{} h$ is defined, then $e \cpsub{} h \in \eqv X$ if and only if $h \in \eqv X$,\label{2_out_of_3_left}
		\item if a pasting $h \subcp{} e$ is defined, then $h \subcp{} e \in \eqv X$ if and only if $h \in \eqv X$.\label{2_out_of_3_right}
	\end{enumerate}
\end{cor}
\begin{proof}
	Follows from Theorem \ref{thm:equivalences_2_out_of_3}, Lemma \ref{lem:weak_inverse_in_Teqv}, and Lemma \ref{lem:composition_in_TA}.
\end{proof}

\begin{prop} \label{prop:all_top_dim_equivalences}
	Let $X$ be a diagrammatic set, $u \in \rd X$, and suppose that every cell $v \submol u$ with $\dim{v} = \dim{u}$ is an equivalence.
	Then $u$ is an equivalence.
\end{prop}
\begin{proof}
	Let $U$ be the shape of $u$, $n \eqdef \dim{U}$, and let $(\order{i}{x})_{i=1}^m$ be an $(n-1)$\nbd ordering induced by an $(n-1)$\nbd layering of $U$.
	For each $i \in \set{1, \ldots, m}$, define $\order{i}{U}$ together with submolecule inclusions $\iota_i\colon \bd{}{-}\order{i}{x} \incl \bd{}{-}\order{i}{U}$ as in the statement of \cite[Proposition 4.3.17]{hadzihasanovic2024combinatorics}, and let $v_i \eqdef \restr{u}{\clset{\order{i}{x}}}$.
	Then
	\[
		\un (\bd{}{-}u) \cp{} u = (\un (\bd{}{-}u) \subcp{\iota_1} v_1) \subcp{\iota_2} v_2 \ldots) \subcp{\iota_m} v_m,
	\]
	and since $\un (\bd{}{-}u) \in \dgn X \subseteq \eqv X$ by Theorem \ref{thm:equivalences_2_out_of_3} and $v_i \in \eqv X$ by assumption, it follows from Corollary \ref{cor:2_out_of_3_consequences}.\ref{2_out_of_3_right} that $\un (\bd{}{-} u) \cp{} u$ is an equivalence and, consequently, that $u$ is an equivalence.
\end{proof}

\begin{dfn}[Equivalent round diagrams]
	Let $u, v$ be a parallel pair of round diagrams in a diagrammatic set $X$.
	We write $u \simeq v$, and say that \emph{$u$ is equivalent to $v$}, if there exists an equivalence $h\colon u \celto v$ in $X$.
\end{dfn}

\begin{prop} \label{prop:equivalence_is_equivalence_relation}
	Let $X$ be a diagrammatic set.
	The relation $\simeq$ is an equivalence relation on $\rd X$.
\end{prop}
\begin{proof}
	Let $u, v, w$ be pairwise parallel round diagrams.
	The unit $\un u\colon u \celto u$ is a degenerate round diagram, so by Theorem \ref{thm:equivalences_2_out_of_3} it exhibits $u \simeq u$.
	Suppose $u \simeq v$, exhibited by an equivalence $e\colon u \celto v$.
	Then by Corollary \ref{cor:2_out_of_3_consequences} there exists an equivalence $e^*\colon v \celto u$, exhibiting $v \simeq u$.
	Finally, suppose $u \simeq v$ and $v \simeq w$, exhibited by $e\colon u \celto v$ and $h\colon v \celto w$, respectively.
	Then by Corollary \ref{cor:2_out_of_3_consequences} $e \cp{} h\colon u \celto w$ is an equivalence, exhibiting $u \simeq w$.
\end{proof}

\begin{prop} \label{prop:equivalence_of_subdiagrams}
	Let $X$ be a diagrammatic set, $u, v, v' \in \rd X$, and $\iota\colon v \submol u$ a rewritable subdiagram.
	If $v \simeq v'$, then $u \simeq \subs{u}{v'}{\iota(v)}$.
\end{prop}
\begin{proof}
	Let $h\colon v \celto v'$ be an equivalence exhibiting $v \simeq v'$.
	Then the pasting $\un {u} \subcp{\iota} h\colon u \celto \subs{u}{v'}{\iota(v)}$ is an equivalence exhibiting $u \simeq \subs{u}{v'}{\iota(v)}$.
\end{proof}

\begin{prop} \label{prop:equivalence_stable_under_equivalence}
	Let $X$ be a diagrammatic set, $u, v \in \rd X$, and suppose $u \simeq v$.
	Then $u$ is an equivalence if and only if $v$ is an equivalence.
\end{prop}
\begin{proof}
	Suppose $v$ is an equivalence and consider an equation $u \cpsub{} x \qeq w$.
	Because $u$ and $v$ are parallel, this also determines an equation $v \cpsub{} x \qeq w$, which admits a lax solution $h\colon v \cpsub{} w' \celto w$ with $h \in \eqv X$.
	Let $k\colon u \celto v$ be an equivalence exhibiting $u \simeq v$ and let $k \cpsub{} h\colon u \cpsub{} w' \celto w$ be the evident pasting.
	This is a lax solution for $u \cpsub{} x \qeq w$, and by Corollary \ref{cor:2_out_of_3_consequences}.\ref{2_out_of_3_left} it is in $\eqv X$.
	The case of an equation $x \subcp{} u \qeq w$ is dual, which proves that $u$ is an equivalence.
	The converse follows by symmetry.
\end{proof}

%%%%%%%%%%%%

\subsection{Weak invertibility and bi-invertibility}

\begin{dfn}[Weakly invertible diagram]
	Let $e\colon u \celto v$ be a round diagram in a diagrammatic set $X$.
	We say that $e$ is \emph{weakly invertible} if there exist a round diagram $e^*\colon v \celto u$ and weakly invertible round diagrams $z\colon e \cp{} e^* \celto \un u$ and $h\colon \un v \celto e^* \cp{} e$.
	In this situation, $e^*$ is called a \emph{weak inverse of $e$}, and $z, h$ are called a \emph{left invertor} and a \emph{right invertor of $e$}, respectively.
\end{dfn}
\noindent We write $\inv X$ for the set of weakly invertible diagrams in $X$.

\begin{comm}
	More formally, $\inv X$ is the greatest fixed point of the operator $\I$ on $\powerset{(\rd X)}$ which sends a set $A \subseteq \rd X$ to the set $\I(A)$ of round diagrams $e\colon u \celto v$ such that there exist round diagrams $e^*\colon v \celto u$, $z\colon e \cp{} e^* \celto \un u$, and $h\colon \un v \celto e^* \cp{} e$ with $z, h \in A$.
	The corresponding coinductive proof method is: if $A \subseteq \I(A)$, then $A \subseteq \inv X$.
\end{comm}

\begin{dfn}[Bi-invertible diagram]
	Let $e\colon u \celto v$ be a round diagram in a diagrammatic set $X$.
	We say that $e$ is \emph{bi-invertible} if there exists a parallel pair of round diagrams $e^L, e^R\colon v \celto u$ and bi-invertible round diagrams $z\colon e \cp{} e^L \celto \un u$ and $h\colon \un v \celto e^R \cp{} e$.
	In this situation, $e^L$ is called a \emph{left inverse} and $e^R$ a \emph{right inverse of $e$}.
\end{dfn}

\noindent We write $\biinv X$ for the set of bi-invertible diagrams in $X$.

\begin{comm}
	The set $\biinv X$ is the greatest fixed point of the operator $\B$ on $\powerset{(\rd X)}$ which sends a set $A \subseteq \rd X$ to the set $\B(A)$ of round diagrams $e\colon u \celto v$ such that there exist round diagrams $e^L, e^R\colon v \celto u$, $z\colon e \cp{} e^L \celto \un u$, and $h\colon \un v \celto e^R \cp{} e$ with $z, h \in A$.
	The corresponding coinductive proof method is: if $A \subseteq \B(A)$, then $A \subseteq \biinv X$.
\end{comm}

\noindent Our next goal is to prove that, in every diagrammatic set $X$,
\[
	\eqv X = \inv X = \biinv X,
\]
that is, all definitions characterise the same class of round diagrams; see \cite{rice2020coinductive, hadzihasanovic2025model} for analogous statements in algebraic models of $\omega$\nbd categories.

\begin{lem} \label{lem:inv_is_biinv}
	Let $X$ be a diagrammatic set, $A \subseteq \rd X$.
	Then $\I(A) \subseteq \B(A)$.
	Consequently, $\inv X \subseteq \biinv X$.
\end{lem}
\begin{proof}
	Let $e \in \rd X$.
	Given data $(e^*, z, h)$ with $z, h \in A$ which exhibit $e \in \I(A)$, the data $(e^L, e^R, z, h)$ with $e^L \equiv e^R \eqdef e^*$ exhibit $e \in \B(A)$.
	It follows that $\inv X = \I(\inv X) \subseteq \B(\inv X)$, so by coinduction $\inv X \subseteq \biinv X$.
\end{proof}

\begin{lem} \label{lem:eqv_is_inv}
	Let $X$ be a diagrammatic set.
	Then $\eqv X \subseteq \inv X$.
\end{lem}
\begin{proof}
	Let $e\colon u \celto v$ be an equivalence in $X$.
	The equations $e \cp{} x \qeq \un u$ and $x \cp{} e \qeq \un v$ admit lax solutions $z\colon e \cp{} e^* \celto \un u$ and $k\colon e' \cp{} e \celto \un v$, where $e^*, e'$ are of type $v \celto u$ and $z, k$ are equivalences.
	Now, we have
	\[
		e^* \simeq \un v \cp{} e^* \simeq e' \cp{} e \cp{} e^* \simeq e' \cp{} \un u \simeq e',
	\]
	with the first and last exhibited by $\lun{}e^*$ and $\run{}e'$, respectively, and the middle ones being instances of Proposition \ref{prop:equivalence_of_subdiagrams}.
	Let $\ell\colon e' \celto e^*$ and $k^*\colon \un v \celto e' \cp{} e$ be equivalences exhibiting $e' \simeq e^*$ and $\un v \simeq e' \cp{} e$, respectively.
	Then the evident pasting $h \eqdef k^* \subcp{} \ell\colon \un v \celto e^* \cp{} e$ is an equivalence by Corollary \ref{cor:2_out_of_3_consequences}.\ref{2_out_of_3_right}.
	We conclude that the data $(e^*, z, h)$ exhibit $e \in \I(\eqv X)$, so by coinduction $\eqv X \subseteq \inv X$.
\end{proof}

\begin{lem} \label{lem:inverses_in_TA}
	Let $X$ be a diagrammatic set, let $e \in \biinv X$, and let $e^L$ be a left inverse of $e$.
	Then $e^L \in \T(\biinv X)$.
\end{lem}
\begin{proof}
	Let $u \eqdef \bd{}{-}e$.
	By definition, there exists a bi-invertible round diagram $z\colon e \cp{} e^L \celto \un u$.
	Then $\un u \in \eqv X \subseteq \biinv X$ by Lemma \ref{lem:inv_is_biinv} and \ref{lem:eqv_is_inv}, and $e, z \in \biinv X$ by assumption.
	It follows that $e^L \in \T(\biinv X)$.
\end{proof}

\begin{lem} \label{lem:biinv_is_eqv}
	Let $X$ be a diagrammatic set.
	Then $\biinv X \subseteq \eqv X$.
\end{lem}
\begin{proof}
	We will prove that $A \eqdef \biinv X$ satisfies the conditions of Lemma \ref{lem:unbiased_main_lemma}.
	First of all, $\dgn X \subseteq \eqv X \subseteq \biinv X$ by Lemma \ref{lem:inv_is_biinv} and \ref{lem:eqv_is_inv}.
	Let $e \in \biinv X$ have type $u \celto v$; by definition, there exist a left and right inverse $e^L, e^R\colon v \celto u$ and bi-invertible $z\colon e \cp{} e^L \celto \un u$ and $h\colon \un v \celto e^R \cp{} e$.
	Let $z^L\colon \un u \celto e \cp{} e^L$ and $h^L\colon e^R \cp{} e \celto \un v$ be left inverses of $z$ and $h$, respectively; by Lemma \ref{lem:inverses_in_TA}, $z^L, h^L \in \T(\biinv X)$.

	Consider an equation $e \cpsub{} x \qeq w$.
	By Remark \ref{rmk:equations_and_subdiagrams}, we have $\iota\colon u \submol \bd{}{-}w$, and we can take the left unitor $\lun{\iota} w\colon w \celto \un u \cpsub{} w$.
	Then the evident pastings
	\[
		k \eqdef z \cpsub{} \rev{(\lun{\iota} w)}\colon e \cpsub{} (e^L \cpsub{\iota} w) \celto w,
		\quad \quad
		k^* \eqdef \lun{\iota} w \subcp{} z^L\colon w \celto e \cpsub{} (e^L \cpsub{\iota} w)
	\]
	exhibit a pair of a lax and colax solution for $e \cpsub{} x \qeq w$, and are both in $\T(\biinv X)$.
	Dually, given an equation $x \subcp{} e \qeq w$, we have $\iota\colon v \submol \bd{}{+}w$, and we can take the right unitor $\run{\iota} w\colon w \subcp{} \un v \celto w$.
	Then the evident pastings
	\[
		k \eqdef h^L \cpsub{} \run{\iota} w\colon (w \subcp{\iota} e^R) \subcp{} e \celto w,
		\quad \quad
		k^* \eqdef \rev{(\run{\iota} w)} \subcp{} h\colon w \celto (w \subcp{\iota} e^R) \subcp{} e
	\]
	exhibit a pair of a lax and colax solution for $x \subcp{} e \qeq w$, and are both in $\T(\biinv X)$.
	This proves that $e \in \E(\T(\biinv X))$, and $\T(\biinv X)$ is unbiased for $e$.
	It follows from Lemma \ref{lem:unbiased_main_lemma} that $\biinv X \subseteq \T(\biinv X) \subseteq \eqv X$.
\end{proof}

\begin{thm} \label{thm:eqv_equals_inv_equals_biinv}
	Let $X$ be a diagrammatic set.
	Then 
	\[ \eqv X = \inv X = \biinv X. \]
\end{thm}
\begin{proof}
	Follows from Lemma \ref{lem:inv_is_biinv}, Lemma \ref{lem:eqv_is_inv}, and Lemma \ref{lem:biinv_is_eqv}.
\end{proof}

\begin{cor} \label{cor:equivalence_is_invertible_up_to_equivalence}
	Let $X$ be a diagrammatic set and $e \in \rd X$ of type $u \celto v$.
	The following are equivalent:
	\begin{enumerate}[label=(\alph*)]
		\item $e$ is an equivalence;
		\item there exists $e^*\colon v \celto u$ such that $e \cp{} e^* \simeq \un u$ and $e^* \cp{} e \simeq \un v$.
	\end{enumerate}
\end{cor}

\begin{rmk}
	The second condition in Corollary \ref{cor:equivalence_is_invertible_up_to_equivalence} is patently symmetric in $e$ and $e^*$, implying that any weak inverse of an equivalence is itself an equivalence.
\end{rmk}

\noindent Let $f\colon X \to Y$ be a morphism of diagrammatic sets and $u\colon U \to X$ a pasting diagram.
We will write $f(u)$ for $f \after u\colon U \to Y$.
Since all of these are defined representably, by precomposition with certain maps in $\rdcpx$, the application $f(-)$ preserves
\begin{enumerate}
	\item boundaries, that is, $f(\bd{n}{\a}u) = \bd{n}{\a}f(u)$,
	\item pastings, that is, $f(u \cpsub{k, \iota} v) = f(u) \cpsub{k, \iota} f(v)$ and $f(v \subcp{k,\iota} u) = f(v) \subcp{k,\iota} f(u)$,
	\item units, that is, $f(\un u) = \un f(u)$,
	\item left and right unitors, that is, $f(\lun{\iota} u) = \lun{\iota} f(u)$ and $f(\run{\iota} u) = \run{\iota} f(u)$,
	\item left and right invertors, that is, $f(\linv u) = \linv f(u)$ and $f(\rinv u) = \rinv f(u)$,
\end{enumerate}
whenever the left-hand side of each equation is defined.

\begin{prop} \label{prop:morphisms_preserve_equivalences}
	Let $f\colon X \to Y$ be a morphism of diagrammatic sets and let $e \in \eqv X$.
	Then $f(e) \in \eqv Y$.
\end{prop}
\begin{proof}
	By Theorem \ref{thm:eqv_equals_inv_equals_biinv}, it suffices to prove that $f(\inv X) \subseteq \I(f(\inv X))$, which by coinduction implies that $f(\inv X) \subseteq \inv Y$.
	Suppose $e$ has type $u \celto v$, and let $e^*$, $z$, $h$ be a weak inverse and invertors of $e$.
	Then $f(z)$ and $f(h)$ have types $f(e) \cp{} f(e^*) \celto \un f(u)$ and $\un f(v) \celto f(e^*) \cp{} f(e)$, respectively, and they belong to $f(\inv X)$.
	We conclude that $f(e) \in \I(f(\inv X))$.
\end{proof}

\begin{cor} \label{cor:morphisms_preserve_simeq}
	Let $f\colon X \to Y$ be a morphism of diagrammatic sets and let $u, v$ be a parallel pair of round diagrams in $X$.
	If $u \simeq v$, then $f(u) \simeq f(v)$.
\end{cor}

%% file: natural.tex
\section{Natural equivalences} \label{sec:natural}

\subsection{Contexts}

\noindent Recall that, if $G$ and $H$ are $\omega$\nbd graphs (in degree $> k$), then a morphism $f\colon G \to H$ is a grade-preserving function that commutes with boundaries.

\begin{dfn}[Context for pasting diagrams]
	Let $X$ be a diagrammatic set.
	For $k$ ranging over $\mathbb{N}$ and $v, w$ over parallel pairs in $\gr{k}{\pd X}$, the class of \emph{contexts on $\pd X(v, w)$} is the inductive class of morphisms of $\omega$\nbd graphs in degree $> k$ with domain $\pd X(v, w)$ generated by the following clauses.
	\begin{enumerate}
		\item (\textit{Left pasting}). 
			For all $u \in \gr{k+1}{\rd X}$ and rewritable $\iota\colon \bd{}{+}u \submol v$,
			\[
				u \cpsub{\iota} -\colon \pd X(v, w) \to \pd X(\subs{v}{\bd{}{-}u}{\iota(\bd{}{+}u)}, w)
			\]
			is a context on $\pd X(v, w)$.
		\item (\textit{Right pasting}).
			For all $u \in \gr{k+1}{\rd X}$ and rewritable $\iota\colon \bd{}{-}u \submol w$,
			\[
				- \subcp{\iota} u\colon \pd X(v, w) \to \pd X(v, \subs{w}{\bd{}{+}u}{\iota(\bd{}{-}u)})
			\]
			is a context on $\pd X(v, w)$.
		\item (\textit{Identity}).
			The identity $-\colon \pd X(v, w) \to \pd X(v, w)$ is a context on $\pd X(v, w)$.
		\item (\textit{Composition}).
			If $\fun{F}\colon \pd X(v, w) \to \pd X(v', w')$ is a context and $\fun{G}$ is a context on $\pd X(v', w')$, then $\fun{GF}$ is a context on $\pd X(v, w)$.
		\item (\textit{Promotion}).
			If $k > 0$ and $\fun{F}$ is a context on $\pd X(\bd{}{-}v, \bd{}{+}w)$, then
			\[
				\fun{F}_{v, w} \eqdef \restr{\fun{F}}{\pd X(v, w)}\colon \pd X(v, w) \to \pd X(\fun{F}v, \fun{F}w)
			\]
			is a context on $\pd X(v, w)$.
	\end{enumerate}
	We let $\dim \fun{F} \eqdef k+1$ be the \emph{dimension} of any context $\fun{F}$ on $\pd X(v, w)$.
\end{dfn}

\begin{rmk}
	It follows from Remark \ref{rmk:strict_omegacats} that, whenever they are well-defined, morphisms of the form $u \cp{n} -$, $u \cpsub{n, \iota} -$, and $u \subcp{n, \iota} -$, as well as their duals and all their iterated composites, are also contexts on $\pd X(v, w)$ for all $n \in \mathbb{N}$ as long as $\dim{u} \leq \dim{v}+1$.
\end{rmk}

\begin{lem} \label{lem:context_layering}
	Let $X$ be a diagrammatic set and $\fun{F}$ a context on $\pd X(v,w)$ with $k \eqdef \dim{\fun{F}}$.
	Then there exist pasting diagrams $(\ell_i, r_i)_{i=1}^{k}$ in $X$ such that
	\begin{enumerate}
		\item $\fun{F} = \ell_k \cp{k-1} (\ell_{k-1} \cp{k-2} \ldots (\ell_1 \cp{0} - \cp{0} r_1) \ldots \cp{k-2} r_{k-1}) \cp{k-1} r_k$,
		\item $\dim{\ell_i}, \dim{r_i} \leq i$ for all $i \in \set{1, \ldots, k}$.
	\end{enumerate}
\end{lem}
\begin{proof}
	We proceed by structural induction on $\fun{F}$.
	If $\fun{F} = u \cpsub{\iota} -$, then the statement is true with
	\[
		\ell_i \eqdef \begin{cases}
			u \cpsub{k-1,\iota} v 
			& \text{if $i = k$}, \\
			\bd{i-1}{-}v
			& \text{if $i < k$},
		\end{cases}
		\quad \quad
		r_i \eqdef \bd{i-1}{+}w
	\]
	for all $i \in \set{1, \ldots, k}$.
	The case of $\fun{F} = - \subcp{\iota} u$ is dual.
	If $\fun{F} = -$, then the statement is true with $\ell_i \eqdef \bd{i-1}{-}v$ and $r_i \eqdef \bd{i-1}{+}w$ for all $i \in \set{1, \ldots, k}$.
	If $\fun{F}$ is obtained by promotion, then it is equal to $\fun{G}_{v, w}$ for some context $\fun{G}$ on $\pd X(\bd{}{-}v, \bd{}{+}w)$ with $\dim{\fun{G}} = k-1$.
	Applying the inductive hypothesis to $\fun{G}$, we obtain $(\ell_i, r_i)_{i=1}^{k-1}$, and the statement is true with $\ell_k \eqdef \fun{G}v$ and $r_k \eqdef \fun{G}w$.
	Finally, suppose $\fun{F}$ is equal to $\fun{HG}$ for some pair of contexts $\fun{G}, \fun{H}$.
	Applying the inductive hypothesis to $\fun{G}$ and $\fun{H}$, we obtain lists of pasting diagrams $(\ell'_i, r'_i)_{i=1}^k$ and $(\ell''_i, r''_i)_{i=1}^k$ in the conditions of the statement.
	Then the statement is true of $\fun{F}$ with $(\ell_i, r_i)_{i=1}^k$ defined by mutual induction by
	\begin{align*}
		\ell_i & \eqdef \begin{cases}
			\ell''_1 \cp{0} \ell'_1 
			& \text{if $i = 1$}, \\
			\ell''_i \cp{i-1} (\ell_{i-1} \cp{i-2} \ell'_i \cp{i-2} r_{i-1})
			& \text{if $i > 1$},
		\end{cases}
		\\
		r_i & \eqdef \begin{cases}
			r'_1 \cp{0} r''_1 
			& \text{if $i = 1$}, \\
			(\ell_{i-1} \cp{i-2} r'_i \cp{i-2} r_{i-1}) \cp{i-1} r''_i
			& \text{if $i > 1$},
		\end{cases}
	\end{align*}
	as can be checked with the axioms of strict $\omega$\nbd categories satisfied by pasting.
\end{proof}

\begin{dfn}[Trim context]
	Let $X$ be a diagrammatic set.
	The class of \emph{trim contexts} is the inductive subclass of contexts generated by left pasting, right pasting, identity, and composition (but not promotion).
\end{dfn}

\begin{lem} \label{lem:trim_context_layering}
	Let $X$ be a diagrammatic set and $\fun{F}$ a context on $\pd X(v, w)$ with $k \eqdef \dim{\fun{F}}$.
	The following are equivalent:
	\begin{enumerate}[label=(\alph*)]
		\item $\fun{F}$ is trim;
		\item there exist pasting diagrams $\ell, r$ in $X$ such that $\fun{F} = \ell \cp{k-1} - \cp{k-1} r$ and $\dim{\ell}, \dim{r} \leq k$.
	\end{enumerate}
\end{lem}
\begin{proof}
	One direction is an easier version of the proof of Lemma \ref{lem:context_layering}.
	For the other direction, suppose $\fun{F} = \ell \cp{k-1} - \cp{k-1} r$, with $\ell$ of shape $L$ and $r$ of shape $R$.
	If $\dim{L} = k$, let $(\order{i}{x})_{i=1}^m$ be a $(k-1)$\nbd ordering of $L$ induced by a $(k-1)$\nbd layering, and let $\order{i}{\ell} \eqdef \restr{\ell}{\clset{\order{i}{x}}}$.
	Then, by \cite[Proposition 4.3.17]{hadzihasanovic2024combinatorics}, 
	\[
		\ell \cp{k-1} - = \order{1}{\ell} \cpsub{} (\order{2}{\ell} \cpsub{} \ldots (\order{m}{\ell} \cpsub{} -)\ldots)
	\]
	which is a trim context on $\pd X(v, w)$.
	If $\dim{L} < k$, then $\ell \cp{k-1} - $ is the identity context, which is evidently trim.
	Proceed dually with $R$.
\end{proof}

\begin{rmk}
	If $\dim \fun{F} = 1$, then $\fun{F}$ is automatically trim.
\end{rmk}

\begin{dfn}[Shape of a context]
	Let $X$ be a diagrammatic set, $v\colon V \to X$ and $w\colon W \to X$ be parallel round diagrams, and $\fun{F}$ be a context on $\pd X(v, w)$ with $k \eqdef \dim{\fun{F}}$.
	Let $(\ell_i\colon L_i \to X, r_i\colon R_i \to X)_{i=1}^k$ be sequences of pasting diagrams provided for $\fun{F}$ by Lemma \ref{lem:context_layering}.
	The \emph{shape of $\fun{F}$} is the molecule 
	\[
		L_k \cp{k-1} (L_{k-1} \cp{k-2} \ldots (L_1 \cp{0} (V \celto W) \cp{0} R_1) \ldots \cp{k-2} R_{k-1}) \cp{k-1} R_k.
	\]
\end{dfn}

\begin{rmk}
	If a cell $a\colon v \celto w$ exists, then the shape of a context $\fun{F}$ on $\pd X(v, w)$ is the shape of $\fun{F}a$.
	Since one can freely attach such a cell to a diagrammatic set when none exist, it follows that the shape of $\fun{F}$ is independent of the choice of $(\ell_i, r_i)_{i=1}^k$.
\end{rmk}

\begin{dfn}[Round context]
	Let $X$ be a diagrammatic set, $v, w$ a parallel pair of round diagrams, and $\fun{F}$ a context on $\pd X(v, w)$.
	We say that $\fun{F}$ is \emph{round} if its shape is round.
\end{dfn}

\begin{rmk}
	If $\fun{F}$ is a round context on $\pd X(v, w)$, then $\fun{F}a$ is round for all round diagrams $a\colon v \celto w$.
\end{rmk}

\begin{rmk} \label{rmk:roundness_and_trimness}
	By \cite[Lemma 7.1.9]{hadzihasanovic2024combinatorics}, if $\fun{F}$ is trim and $v, w$ are round, then $\fun{F}$ is round.
\end{rmk}

\begin{dfn}[Weakly invertible context]
	Let $X$ be a diagrammatic set.
	The class of \emph{weakly invertible contexts} is the inductive subclass of contexts obtained by restricting left pasting and right pasting to $u \in \eqv X$.
\end{dfn}

\begin{lem} \label{lem:context_factorisation}
	Let $X$ be a diagrammatic set, let $\fun{F}$ be a context on $\pd X(v, w)$, and suppose $\dim \fun{F} > 1$.
	Then there exist
	\begin{enumerate}
		\item a context $\fun{G}$ on $\pd X(\bd{}{-}v, \bd{}{+}w)$, and
		\item a trim context $\fun{T}$ on $\pd X(\fun{G}v, \fun{G}w)$
	\end{enumerate}
	such that $\fun{F} = \fun{T}\fun{G}_{v, w}$.
	Moreover,
	\begin{enumerate}
		\item if $\fun{F}$ is round, then $\fun{T}$ and $\fun{G}$ are round,
		\item if $\fun{F}$ is weakly invertible, then $\fun{T}$ and $\fun{G}$ can be chosen weakly invertible.
	\end{enumerate}
\end{lem}
\begin{proof}
	The main statement follows from Lemma \ref{lem:context_layering} and Lemma \ref{lem:trim_context_layering} with
	\begin{enumerate}
		\item $\fun{G} \eqdef \ell_{k-1} \cp{k-2} \ldots (\ell_1 \cp{0} - \cp{0} r_1) \ldots \cp{k-2} r_{k-1}$,
		\item $\fun{T} \eqdef \ell_k \cp{k-1} - \cp{k-1} r_k$.
	\end{enumerate}
	Suppose that $\fun{F}$ is round, and let
	\[
		U' \eqdef L_{k-1} \cp{k-2} \ldots (L_1 \cp{0} (V \celto W) \cp{0} R_1) \ldots \cp{k-2} R_{k-1}
	\]
	so that the shape of $\fun{F}$ is $U \eqdef L_k \cp{k-1} U' \cp{k-1} R_k$.
	Then the shape of $\fun{G}$ is $\subs{\bd{}{-}U'}{\compos{V}}{V}$, which is round if and only if $\bd{}{-}U'$ is round.
	If $\dim{L_k} < k$, then $\bd{}{-}U' = \bd{}{-}U$, which is round by assumption.
	Otherwise, we can take a $(k-1)$\nbd ordering of $L_k$ induced by a $(k-1)$\nbd layering, and use \cite[Corollary 4.3.15]{hadzihasanovic2024combinatorics} to deduce that $\bd{}{-}U'$ is round from the fact that $\bd{}{-}U = \bd{}{-}L_k$ is round.
	This proves that $\fun{G}$ is round.
	Since $\fun{G}v$ and $\fun{G}w$ are round and $\fun{T}$ is trim, it follows from Remark \ref{rmk:equations_and_subdiagrams} that $\fun{T}$ is also round.

	Finally, suppose that $\fun{F}$ is weakly invertible; we proceed by structural induction.
	If $\fun{F}$ is produced by left pasting, right pasting, or identity, we can take $\fun{T} \eqdef \fun{F}$ and $\fun{G} \eqdef -$, both weakly invertible.
	If $\fun{F}$ is produced by promotion, it is equal to $\fun{G}_{v, w}$ for some weakly invertible $\fun{G}$, and we can take $\fun{T} \eqdef -$.
	If $\fun{F}$ is produced by composition, then by the inductive hypothesis it is equal to $\fun{T'}\fun{G'}_{v',w'}\fun{T}\fun{G}_{v,w}$ with $\fun{T}, \fun{T'}, \fun{G}, \fun{G'}$ weakly invertible and $\fun{T}, \fun{T'}$ trim.
	Then a similar argument as in the proof of Lemma \ref{lem:context_layering}, using distributivity of lower-dimensional pasting over higher-dimensional pasting, proves that we can rewrite $\fun{G'}_{v',w'}\fun{T}$ as $\fun{T''}\fun{G''}_{v'', w''}$, where both $\fun{T''}$ and $\fun{G''}$ are weakly invertible and $\fun{T''}$ is trim.
	This completes the proof.
\end{proof}

%%%%%%%%%%%%

\subsection{Natural equivalences of round contexts}

\begin{dfn}[Natural equivalence of round contexts]
	Let $X$ be a diagrammatic set and let $\fun{F}, \fun{G}\colon \pd X(v, w) \to \pd X(v', w')$ be round contexts.
	A family of equivalences $\th a\colon \fun{F}a \celto \fun{G}a$ indexed by round diagrams $a\colon v \celto w$ is a \emph{natural equivalence from $\fun{F}$ to $\fun{G}$} if, for all round diagrams $a, b\colon v \celto w$, there exists a natural equivalence from ${\fun{F}_{a,b}-} \cp{} \th b$ to $\th a \cp{} {\fun{G}_{a,b}-}$ as round contexts $\pd X(a, b) \to \pd X(\fun{F}a, \fun{G}b)$.
	We write $\th\colon \fun{F} \eqvto \fun{G}$ to indicate that $\th$ is a natural equivalence from $\fun{F}$ to $\fun{G}$.
\end{dfn}

\begin{comm}
	This is another coinductive definition, only slightly complicated by the fact that it determines not a set of diagrams, but a set of families of diagrams for each parallel pair of round contexts.
	Nevertheless, the principle is the same.
	Let $\ctxeqv X$ denote the set of triples $(\fun{F}, \fun{G}, \th)$ such that
	\begin{enumerate}
		\item $\fun{F}, \fun{G}$ are parallel round contexts $\pd X(v, w) \to \pd X(v', w')$,
		\item $\th$ is a family of equivalences $\th a\colon \fun{F}a \celto \fun{G}a$ indexed by round diagrams $a\colon v \celto w$.
	\end{enumerate}
	Then, for each subset $A \subseteq \ctxeqv X$, let $\N(A)$ be the set of triples $(\fun{F}, \fun{G}, \th)$ such that, for all indexing round diagrams $a, b$ of $\th$, there exists a triple $({\fun{F}_{a,b}-} \cp{} \th b, \th a \cp{} {\fun{G}_{a,b}-}, \psi )$ in $A$.
	This determines an order-preserving operator $\N$ on $\powerset{(\ctxeqv X)}$, whose greatest fixed point is, by definition, the set $\nateqv X$ of natural equivalences.
\end{comm}

\noindent Our main goal in this section will be to prove that natural equivalences contain certain special families of equivalences, and are closed under a number of operations.
We start by introducing these special families: the \emph{higher unitors} and the \emph{context pushforwards}.

\begin{dfn}[Higher unitors]
	Let $u: U \to X$ be a pasting diagram in a diagrammatic set, $k < \dim{u}$, and let $\iota\colon v \submol \bd{k}{-}u$ and $j\colon w \submol \bd{k}{+}u$ be rewritable subdiagrams of shapes $V$ and $W$, respectively.
	Let $K \eqdef \bd{k}{}U \setminus \inter{\iota(V)}$ and $K' \eqdef \bd{k}{}U \setminus \inter{j(W)}$.
	We define
	\begin{itemize}
		\item $\lun{k, \iota}u$ to be the degenerate pasting diagram $u \after \tau_K\colon \arr \pcyl{K} U \to X$,
		\item $\run{k, j}u$ to be the degenerate pasting diagram $u \after \tau_{K'}\colon \arr \pcyl{K'} U \to X$.
	\end{itemize}
	The types of $\lun{k, \iota}u$ and $\run{k, j}u$ are defined inductively as follows, together with subdiagrams $u \submol \bd{}{\a} (\lun{k, \iota}u)$ and $u \submol \bd{}{\a} (\run{k, j}u)$ for all $\a \in \set{+, -}$:
	\begin{itemize}
		\item if $k = \dim{u} - 1$, then 
			\[
				\lun{k, \iota}u \equiv \lun{\iota} u\colon u \celto \un v \cpsub{\iota} u, \quad \quad
				\run{k, j} u \equiv \run{j} u\colon u \subcp{j} \un w \celto u, 
			\]
			the previously defined left and right unitors of $u$;
		\item if $k < \dim{u} - 1$, then 
			\[
				\lun{k, \iota}u\colon u \cpsub{} \lun{k, \iota}(\bd{}{+}u) \celto \lun{k, \iota}(\bd{}{-}u) \subcp{} u, \;
				\run{k, j}u\colon u \cpsub{} \run{k, j}(\bd{}{+}u) \celto \run{k, j}(\bd{}{-}u) \subcp{} u,
			\]
			where the pastings are at the specified subdiagrams $\bd{}{\a}u \submol \bd{}{\b}\lun{k, \iota}(\bd{}{\a}u)$ and $\bd{}{\a}u \submol \bd{}{\b}\run{k, j}(\bd{}{\a}u)$ for all $\a, \b \in \set{ +, - }$.
	\end{itemize}
	The subdiagram inclusions are the evident ones at each step.
\end{dfn}

\noindent
We let $\unitor X \subseteq \ctxeqv X$ be the set of triples $(\fun{F}, \fun{G}, \th)$ such that, for some $n \in \mathbb{N}$ and parallel $v, w \in \gr{n}{\rd X}$, both $\fun{F}$ and $\fun{G}$ are round contexts on $\pd X(v, w)$ with
\begin{itemize}
	\item $\fun{F} = \fun{G} = -$ and $\th a = \un a$, or
	\item there exists a rewritable $\iota\colon u \submol v$ such that $\fun{F} = -$, $\fun{G} = \un{u} \cpsub{\iota} -$, and $\th a = \lun{\iota} a$, or
	\item there exists a rewritable $j\colon u \submol w$ such that $\fun{F} = - \subcp{j} \un{u}$, $\fun{G} = -$, and $\th a = \run{j} a$, or
	\item there exists $k < n$ and a rewritable $\iota\colon u \submol \bd{k}{-}v = \bd{k}{-}w$ such that $\fun{F} = - \cpsub{} \lun{k,\iota}w$, $\fun{G} = \lun{k,\iota}v \subcp{} -$, and $\th a = \lun{k, \iota} a$, or
	\item there exists $k < n$ and a rewritable $j\colon u \submol \bd{k}{+}v = \bd{k}{+}w$ such that $\fun{F} = - \cpsub{} \run{k,j}w$, $\fun{G} = \run{k,j}v \subcp{} -$, and $\th a = \run{k, \iota} a$,
\end{itemize}
where $a$ ranges over round diagrams of type $v \celto w$.

\begin{exm}
	To motivate the definition of higher unitors, we demonstrate informally how they will be used to show that contexts of the form \( \un x \cpsub{} - \) are naturally equivalent to identity contexts.
	In this example, \( \iota \) will always denote the subdiagram \( \iota \colon \bd{0}{-}w \submol w \) for a pasting diagram \( w \) that the reader can infer from the context.
	Suppose that \( x \) is a \( 0 \)\nbd cell in a diagrammatic set \( X \), and let \( \fun{F} \) be the context \( \un x \cpsub{} - \) on \( \rd X(x, y) \) for a \( 0 \)\nbd cell \( y \). 
	Then \( \fun{F} \) is parallel to the identity context \( \fun{I} \colon \rd X(x, y) \to \rd X(x, y) \).
	The family of degenerate diagrams \( \lun{0, \iota} a \colon a \celto \un x \cpsub{} a \), indexed by round diagrams \( a \colon x \celto y \), is a family of equivalences of type \( \lun{0, \iota} a \colon \fun{I} a \celto \fun{F}a \).
	Since we want this family to be natural, for all round diagrams \( a, b \colon x \celto y \), we need a family of equivalences of type \( \fun{I}_{a, b} w \cp{} \lun{0, \iota} b \celto \lun{0, \iota}a \cp{} \fun{F}_{a, b}w \), indexed by round diagrams \( u \colon a \celto b \).
	These equivalences are given by the next higher unitors \( \lun{0, \iota} u \), which are \( 3 \)\nbd dimensional round diagrams whose input and output \( 2 \)\nbd boundaries are
	\begin{center}
		\begin{tikzcd}[row sep=large]
			& x &&&& x \\
			x && y & {\text{and}} & x && {y.}
			\arrow["b"{description}, curve={height=-6pt}, from=1-2, to=2-3]
			\arrow[""{name=0, anchor=center, inner sep=0}, "b"{description}, curve={height=-12pt}, from=1-6, to=2-7]
			\arrow[""{name=1, anchor=center, inner sep=0}, "a"{description}, curve={height=12pt}, from=1-6, to=2-7]
			\arrow["{\un x}"{description}, curve={height=-6pt}, from=2-1, to=1-2]
			\arrow[""{name=2, anchor=center, inner sep=0}, "b"{description}, curve={height=-6pt}, from=2-1, to=2-3]
			\arrow[""{name=3, anchor=center, inner sep=0}, "a"{description}, curve={height=18pt}, from=2-1, to=2-3]
			\arrow["{\un x}"{description}, curve={height=-6pt}, from=2-5, to=1-6]
			\arrow[""{name=4, anchor=center, inner sep=0}, "a"{description}, curve={height=18pt}, from=2-5, to=2-7]
			\arrow["u", between={0.2}{0.8}, Rightarrow, from=1, to=0]
			\arrow["{\lun{0, \iota}b}"{description}, between={0.2}{1}, Rightarrow, from=2, to=1-2]
			\arrow["u", between={0.2}{0.8}, Rightarrow, from=3, to=2]
			\arrow["{\lun{0, \iota}a}"{description, pos=0.4}, between={0.2}{0.7}, Rightarrow, from=4, to=1-6]
		\end{tikzcd}
	\end{center}
	Iterating this construction, one concludes that the contexts \( \fun{I} \) and \( \fun{F} \) are naturally equivalent.
\end{exm}

\begin{dfn}[Context subdiagram]
	Let $X$ be a diagrammatic set, $v, w \in \pd X$ be parallel, and let $\fun{F}$ be a context on $\pd X(v, w)$.
	A \emph{context subdiagram} $\iota\colon z \submol \fun{F}$ is a pair of 
	\begin{enumerate}
		\item a decomposition $\fun{F} = v' \cp{} {\fun{F}'-}$ or $\fun{F} = {\fun{F}'-} \cp{} v'$, and
		\item a subdiagram $\iota\colon z \submol v'$.
	\end{enumerate}
	A context subdiagram is \emph{rewritable} if $\dim{v'} = \dim \fun{F}$ and $\iota$ is rewritable.
\end{dfn}

\begin{lem} \label{lem:context_subdiagrams_give_subdiagrams}
	Let $X$ be a diagrammatic set, $v, w \in \pd X$ be parallel, let $\fun{F}$ be a context on $\pd X(v, w)$, and let $\iota\colon z \submol \fun{F}$ be a context subdiagram.
	Then, for all $a\colon v \celto w$, $\iota$ determines a subdiagram $\iota_a\colon z \submol \fun{F}a$, which is rewritable if $\iota$ is rewritable.
\end{lem}
\begin{proof}
	Let \( \iota \colon z \submol \fun{F} \) be a context subdiagram. 
	By definition, \( \fun{F} = v' \cp{} {\fun{F}'-} \) or \( \fun{F} = {\fun{F}'-} \cp{} v' \) with $\iota\colon z \submol v'$; without loss of generality, suppose we are in the first case.
	Then, for all \( a \colon v \celto w \), \( {\fun{F}a} = v' \cp{} {\fun{F'}a} \), and composing \( \iota \) with the evident subdiagram $v' \submol \fun{F}a$ determines a subdiagram of \( \iota_a\colon z \submol {\fun{F}a} \), which is rewritable if \( \iota \) is.
\end{proof}

\begin{dfn}[Context pushforward]
	Let $X$ be a diagrammatic set, $v, w \in \rd X$ be parallel, $\fun{F}$ be a round context on $\pd X(v, w)$, and let $\iota\colon z \submol \fun{F}$ be rewritable.
	Given an equivalence $h\colon z \celto z'$, the \emph{context pushforward of $\fun{F}$ along $h$ at $\iota$} is the family $\un{(\fun{F})} \subcp{\iota} h$ of equivalences
	\[
		\un{(\fun{F}a)} \subcp{\iota_a} h\colon \fun{F}a \celto \subs{\fun{F}a}{z'}{\iota_a(z)}
	\]
	indexed by round diagrams $a\colon v \celto w$.
	We let $\ctxpfw X \subseteq \ctxeqv X$ be the set of triples
\[
	(\fun{F}, \subs{{\fun{F}}}{z'}{\iota(z)}, \un{(\fun{F})} \subcp{\iota} h)
\]
for some round context $\fun{F}$, rewritable context subdiagram $\iota\colon z \submol \fun{F}$, and equivalence $h\colon z \celto z'$.
\end{dfn}

\begin{exm}
	Let \( z \colon u \celto v \) be a \( 1 \)\nbd cell in a diagrammatic set, let \( \fun{F} \) be the context \( z \cpsub{} - \) on \( \pd X(v, w) \) for a \( 0 \)\nbd cell \( w \), and \( h \colon z \celto z' \) be a \( 2 \)\nbd cell. 
	Then \( \fun{F} \) has an evident rewritable context subdiagram \( \iota \colon z \submol \fun{F} \), and, substituting \( z' \) for \( z \), we have the context \( \subs{{\fun{F}}}{z'}{\iota(z)} \):
	\begin{center}
		\begin{tikzcd}
			{\subs{{\fun{F}}}{z'}{\iota(z)} } & u & v & w \\
			{\fun{F}} & u & v & w.
			\arrow["{z'}"{description}, from=1-2, to=1-3]
			\arrow[dotted, from=1-3, to=1-4]
			\arrow["z"{description}, from=2-2, to=2-3]
			\arrow[dotted, from=2-3, to=2-4]
		\end{tikzcd}
	\end{center}
	If \( h \) is an equivalence, then for each \( a \colon v \celto w \), the round diagram
	\begin{equation*}
		\un{(\fun{F}a)} \subcp{\iota_a} h\colon \fun{F}a \celto \subs{\fun{F}a}{z'}{\iota_a(z)}
	\end{equation*}
	\begin{center}
		\begin{tikzcd}
			u && v & w \\
			\\
			& v
			\arrow[""{name=0, anchor=center, inner sep=0}, "z"{description}, curve={height=12pt}, from=1-1, to=1-3]
			\arrow[""{name=1, anchor=center, inner sep=0}, "{z'}"{description}, curve={height=-12pt}, from=1-1, to=1-3]
			\arrow["z"{description}, curve={height=18pt}, from=1-1, to=3-2]
			\arrow["a", from=1-3, to=1-4]
			\arrow["{\un{v}}"{description}, from=3-2, to=1-3]
			\arrow[""{name=2, anchor=center, inner sep=0}, "a"{description}, curve={height=18pt}, from=3-2, to=1-4]
			\arrow["h", between={0.2}{0.8}, Rightarrow, from=0, to=1]
			\arrow["{\lun{}a}"{description}, between={0.2}{1}, Rightarrow, from=2, to=1-3]
			\arrow["{\run{}a}"{description}, between={0}{0.8}, Rightarrow, from=3-2, to=0]
		\end{tikzcd}
	\end{center}
	is also an equivalence.
	This family of equivalences \( \un(\fun{F}) \subcp{\iota} h \) is the context pushforward of \( \fun{F} \) along \( h \) at \( \iota \).
\end{exm}

\noindent Next, given $A \subseteq \ctxeqv X$, we let $\T(A)$ denote the closure of $A$ under the following clauses.
\begin{enumerate}
	\item (\textit{Composition}). If $(\fun{F}, \fun{G}, \th), (\fun{G}, \fun{H}, \psi) \in \T(A)$, then $(\fun{F}, \fun{H}, \th \cp{} \psi) \in \T (A)$, where $\th \cp{} \psi{}$ is defined by
		\[
			(\th \cp{} \psi{})a \eqdef \th a \cp {} \psi a\colon \fun{F} a \celto \fun{H} a
		\]
		for all indexing round diagrams $a$.
	\item (\textit{Left context}). If $(\fun{F}, \fun{G}, \th) \in \T(A)$, with $\fun{F}, \fun{G}\colon \pd X(v, w) \to \pd X(v', w')$ parallel round contexts, and $\fun{H}$ is a round context on $\pd X(v', w')$, then $(\fun{HF}, \fun{HG}, \fun{H}\th \cp{} \un ({\fun{HG}})) \in \T(A)$, where the family $\fun{H}\th \cp{} \un({\fun{HG}})$ assigns to each round diagram $a\colon v \celto w$ the equivalence
		\[
			\fun{H} \th a \cp{} \un{(\fun{HG}a)}\colon \fun{HF}a \celto \fun{HG}a.
		\]
	\item (\textit{Right context}). If $(\fun{F}, \fun{G}, \th) \in \T(A)$, with $\fun{F}, \fun{G}$ parallel round contexts on $\pd X(v', w')$, and $\fun{H}\colon \pd X(v, w) \to \pd X(v', w')$ is a round context, then $(\fun{FH}, \fun{GH}, \th{\fun{H}}) \in \T(A)$, where the family $\th\fun{H}$ assigns to each round diagram $a\colon v \celto w$ the equivalence
		\[
			\th \fun{H}a\colon \fun{FH}a \celto \fun{GH}a.
		\]
	\item (\textit{Weak inversion}). If $(\fun{F}, \fun{G}, \th) \in \T(A)$, then $(\fun{G}, \fun{F}, \th^*) \in \T(A)$ for every choice of componentwise weak inverses
		\[
			\th^*a\colon \fun{G} a \celto \fun{F} a
		\]
		to $\th a$ at the indexing round diagrams $a$.
\end{enumerate}
Our next goal is to prove that 
\[
	\nateqv X = \T(\nateqv X \cup \unitor X \cup \ctxpfw X),
\]
that is, higher unitors and context pushforwards are natural equivalences, and natural equivalences are closed under composition, context, and weak inversion.
We will follow a similar strategy to the one we used to prove saturation properties of equivalences.

\begin{lem} \label{lem:unitors_and_pushforwards_in_NT}
	Let $X$ be a diagrammatic set and $A \eqdef \unitor X \cup \ctxpfw X$.
	Then $A \subseteq \N(\T(A))$.
\end{lem}
\begin{proof}
	Let $(\fun{F}, \fun{G}, \th) \in A$, where $\fun{F}$ and $\fun{G}$ are round contexts on $\pd X(v, w)$, and let $n \eqdef \dim{v}$.
	Suppose \( \fun{F} = \fun{G} = - \) and \( \th a = \un a \).
	Then \( (- \cp{} \un b, -, \run{}) \) and \( (-, \un a \cp{} -, \lun{}) \) are both in \( A \), so by closure under composition the triple \( (- \cp{} \un b, \un a \cp{} -, \run{} \cp{} \lun{}) \) is in \( \T(A) \).
	Next, if \( \fun{F} = - \), \( \fun{G} = \un u \cpsub{\iota} - \), and \( \th = \lun{\iota} \), then for all round diagrams $a, b\colon v \celto w$, 
	\[
		{\fun{F}_{a,b}-} \cp{} \th b = - \cp{} \lun{\iota}b,
		\quad \quad
		\th a \cp{} {\fun{G}_{a,b}-} = \lun{\iota}a \cp{} -,
	\]
	and $(- \cp{} \lun{\iota}b, \lun{\iota}a \cp{}-, \lun{n, \iota}) \in A$.
	If $\fun{F} = - \cpsub{} \lun{k, \iota}w$, $\fun{G} = \lun{k, \iota}v \subcp{} -$ and $\th = \lun{k, \iota}$ for some $k < n$, then for all round diagrams $a, b\colon v \celto w$,
	\begin{align*}
		{\fun{F}_{a,b}-} \cp{} \th b & = ({-} \cpsub{} \lun{k, \iota}w) \cp{} \lun{k, \iota}b = {-} \cpsub{} \lun{k, \iota}b, \\
		\th a \cp{} {\fun{G}_{a,b}-} & = \lun{k, \iota}a \cp{} (\lun{k, \iota}w \subcp{} {-}) = \lun{k, \iota}a \subcp{} {-},
	\end{align*}
	and $({-} \cpsub{} \lun{k, \iota}b, \lun{k,\iota}a \subcp{} -, \lun{k,\iota}) \in A$.
	The case of right unitors is dual.

	Finally, suppose $\fun{G} = \subs{\fun{F}}{z'}{\iota(z)}$ and $\th = \un{(\fun{F})} \subcp{\iota} h$ for some rewritable context subdiagram $\iota\colon z \submol \fun{F}$ and equivalence $h\colon z \celto z'$.
	Then, for all round diagrams $a, b\colon v \celto w$,
	\begin{align*}
		{\fun{F}_{a,b}-} \cp{} \th b & = {\fun{F}_{a,b}-} \cp{} (\un{\fun{F}b} \subcp{\iota_b} h) = ({\fun{F}_{a,b}-} \cp{} \un{\fun{F}b}) \subcp{\iota_b} h, 
		\\
	\th a \cp{} {\fun{G}_{a,b}-} & = (\un{\fun{F}a} \subcp{\iota_a} h) \cp{} \subs{{\fun{F}_{a,b}-}}{z'}{\iota(z)} = (\un{\fun{F}a} \cp{} {\fun{F}_{a,b}-}) \subcp{\iota_b} h.
	\end{align*}
	Now, by the first part of the proof and by closure under the right context $\fun{F}_{a,b}$,
	\[
		({\fun{F}_{a,b}-} \cp{} \un{\fun{F}b}, \un{\fun{F}a} \cp{} {\fun{F}_{a,b}-}, (\run{} \cp{} \lun{})\fun{F}_{a,b}) \in \T(A),
	\]
	and we conclude by closure under the left context $- \subcp{\iota_b} h$.
\end{proof}

\begin{lem} \label{lem:natural_saturation_main_lemma}
	Let $X$ be a diagrammatic set and let $A \subseteq \ctxeqv X$ such that $\unitor X \cup \ctxpfw X \subseteq A$.
	If $A \subseteq \N(\T(A))$, then $\T(A) \subseteq \nateqv X$.
\end{lem}
\begin{proof}
	We have \( \T(A) \subseteq \T(\N(\T(A))) \), so to conclude by coinduction it suffices to show that $ \T(\N(\T(A))) = \N(\T(A))$.
	
	We introduce the following relation: if \( \fun{F}, \fun{G} \) are round contexts, we write \( \fun{F} \sim \fun{G} \) if and only if there exists \( \th \) such that \( (\fun{F}, \fun{G}, \th) \in \T(A) \). 
	We claim that \( \sim \) is a congruence on round contexts with respect to composition.
	By Lemma \ref{lem:unitors_and_pushforwards_in_NT}, $(-, -, \un{}) \in \T(A)$, so given any round context $\fun{F}$, closure under left context implies that $(\fun{F}, \fun{F}, \un{\fun{F}}) \in \T(A)$, that is, $\fun{F} \sim \fun{F}$.
	If $\fun{F} \sim \fun{G}$, then closure of $\T(A)$ under weak inversion implies $\fun{G} \sim \fun{F}$.
	If $\fun{F} \sim \fun{G}$ and $\fun{G} \sim \fun{H}$, then closure of $\T(A)$ under composition implies $\fun{F} \sim \fun{H}$, which proves that $\sim$ is an equivalence relation.
	Finally, by closure under left and right context, $\fun{F} \sim \fun{G}$ implies $\fun{HFK} \sim \fun{HGK}$ for any pair of round contexts $\fun{H}, \fun{K}$ that can be respectively post-composed and pre-composed with $\fun{F}$ and $\fun{G}$.

	Now, observe that, if $\fun{F}$ and $\fun{G}$ are round contexts on $\pd X(v, w)$, we have $(\fun{F}, \fun{G}, \th) \in \N(\T(A))$ if and only if, for all round diagrams $a, b\colon v \celto w$, we have ${\fun{F}_{a,b}-} \cp{} \th b \sim \th a \cp{} {\fun{G}_{a,b}-}$.
	Suppose that \( (\fun{F}, \fun{G}, \th), (\fun{G}, \fun{H}, \psi) \in \N(\T(A)) \), where $\fun{F}$ is a round context on $\pd X(v,w)$.  
	Then, for all round diagrams \( a, b\colon v \celto w \),
	\[
		{\fun{F}_{a,b}-} \cp{} \th b \cp{} \psi b 
		\sim \th a \cp{} {\fun{G}_{a,b}-} \cp{} \psi b
		\sim \th a \cp{} \psi a \cp{} {\fun{H}_{a,b}-},
	\]
	which proves that $\N(\T(A))$ is closed under composition.
	Next, suppose \( (\fun{F}, \fun{G}, \th) \in \N(\T(A)) \), where $\fun{F}, \fun{G}\colon \pd X(v, w) \to \pd X(v', w')$ are round contexts, and let \( \fun{H} \) be a round context on $\pd X(v', w')$.
	For all round diagrams $a, b\colon v \celto w$, we have
	\begin{align*}
		{(\fun{HF})_{a,b}-} \cp{} \fun{H}\th b \cp{} \un (\fun{H}\fun{G} b) 
		& = \fun{H}_{\fun{F}a,\fun{G}b}({\fun{F}_{a,b}-} \cp{} \th b) \cp{} \un(\fun{H}\fun{G}b) \\
		& \sim \fun{H}_{\fun{F}a,\fun{G}b}(\th a \cp{} {\fun{G}_{a,b}-}) \cp{} \un(\fun{H}\fun{G}b) \\
		& = \fun{H}\th a \cp{} {(\fun{HG})_{a,b}-} \cp{} \un(\fun{H}\fun{G}b) \\
		& \sim \fun{H}\th a \cp{} \un(\fun{HG} a) \cp{} {(\fun{HG})_{a,b}-} \cp{} \un(\fun{H}\fun{G}b) \\
		&\sim \fun{H}\th a \cp{} \un (\fun{HG} a) \cp{} (\fun{HG})_{a,b}- 
	\end{align*}
	where the identities are proved using Lemma \ref{lem:context_layering} on $\fun{H}$, while the final two steps are applications of a left and a right unitor in context, respectively.
	This proves that $\N(\T(A))$ is closed under left context; an analogous argument proves that it is closed under right context.

	Finally, suppose \( (\fun{F}, \fun{G}, \th) \in \N(\T(A)) \), where $\fun{F}$ and $\fun{G}$ are round contexts on $\pd X(v, w)$, and let \( \th^* \) be a choice of componentwise weak inverses for \( \th \).
	For each round diagram $a\colon v \celto w$, let $z_a\colon \th a \cp{} \th^* a \celto \un{\fun{F}a}$, $h_a\colon \un{\fun{G}a} \celto \th^*a \cp{} \th a$ be a choice of left invertor and right invertor for $\th a$.
	Then, for all round diagrams $a, b\colon v \celto w$, we have
	\begin{align*}
		{\fun{G}_{a,b}-} \cp{} \th^* b 
		&\sim \un(\fun{G}a) \cp{} {\fun{G}_{a,b}-} \cp{} \th^* b 
			&& \text{by left unitor} \\
		&\sim \th^* a \cp{} \th a \cp{} {\fun{G}_{a,b}-} \cp{} \th^* b 
			&& \text{by pushforward with $h_a$} \\
		&\sim \th^* a \cp{} {\fun{F}_{a,b}-} \cp{} \th b \cp{} \th^* b \\
		&\sim \th^* a \cp{} {\fun{F}_{a,b}-} \cp{} \un{\fun{F} b} 
			&& \text{by pushforward with $z_b$} \\
		&\sim \th^* a \cp{} {\fun{F}_{a,b}-} 
			&& \text{by right unitor},
	\end{align*}
	since $\T(A)$ contains context pushforwards and unitors.
	This proves that $\N(\T(A))$ is closed under weak inversion, which completes the proof.
\end{proof}

\begin{thm} \label{thm:natural_saturation}
	Let $X$ be a diagrammatic set.
	Then
	\begin{enumerate}
		\item $\unitor X \cup \ctxpfw X \subseteq \nateqv X$, and
		\item $\nateqv X = \T(\nateqv X)$.
	\end{enumerate}
\end{thm}
\begin{proof}
	Let $A \eqdef \nateqv X \cup \unitor X \cup \ctxpfw X$.
	Then $A \subseteq \N(\T(A))$ by the fact that $\nateqv X = \N(\nateqv X) \subseteq \N(\T(\nateqv X))$ combined with Lemma \ref{lem:unitors_and_pushforwards_in_NT}.
	We conclude by Lemma \ref{lem:natural_saturation_main_lemma}.
\end{proof}

\begin{dfn}[Equivalent round contexts]
	Let $X$ be a diagrammatic set and let $\fun{F}, \fun{G}$ be round contexts $\pd X(v, w) \to \pd X(v', w')$.
	We write $\fun{F} \simeq \fun{G}$, and say that \emph{$\fun{F}$ is equivalent to $\fun{G}$}, if there exists a natural equivalence $\th\colon \fun{F} \eqvto \fun{G}$.
\end{dfn}

\begin{prop} \label{prop:natural_equivalence_is_equivalence_relation}
	Let $X$ be a diagrammatic set. 
	Then the relation $\simeq$ on round contexts in $X$ is
	\begin{enumerate}
		\item an equivalence relation,
		\item a congruence with respect to composition of round contexts,
		\item compatible with the relation $\simeq$ on round diagrams of the same dimension, that is, if $v \simeq w$ and $\fun{F} \simeq \fun{G}$ with $\dim{v} = \dim{\fun{F}}$, then
			\begin{itemize}
				\item if $v \cp{} {\fun{F}-}$ is defined, then $v \cp{} {\fun{F}-} \simeq w \cp{} {\fun{G}-}$,
				\item if ${\fun{F}-} \cp{} v$ is defined, then ${\fun{F}-} \cp{} v \simeq {\fun{G}-} \cp{} w$.
			\end{itemize}
	\end{enumerate}
\end{prop}
\begin{proof}
	The proof that $\simeq$ is a congruence is the same as the proof that $\sim$ is a congruence in Lemma \ref{lem:natural_saturation_main_lemma}, specialised to $A \eqdef \nateqv X$, which is admissible by Theorem \ref{thm:natural_saturation}.
	Suppose that $v \simeq w$ and $\fun{F} \simeq \fun{G}$, where $v, w$ are round diagrams and $\fun{F}, \fun{G}$ round contexts with $\dim{\fun{F}} = \dim{v}$, and suppose $v \cp{} {\fun{F}-}$ is defined.
	Let $h\colon v \celto w$ be an equivalence exhibiting $v \simeq w$.
	Then
	\begin{align*}
		v \cp{} {\fun{F}-}
		&\simeq w \cp{} {\fun{F}-}
			&& \text{by context pushforward with $h$} \\
		&\simeq w \cp{} {\fun{G}-}
	\end{align*}
	since context pushforwards are natural equivalences by Theorem \ref{thm:natural_saturation}.
	The case where ${\fun{F}-} \cp{} v$ is defined is dual.
\end{proof}

%%%%%%%%

%% file: bicategory.tex
\section{Bicategories of round diagrams} \label{sec:bicategory}

Let $X$ be a diagrammatic set.
We will show that, for each $n \in \mathbb{N}$, one can form a strictly associative bicategory $\bic{n}{X}$ whose set of 0\nbd cells is $\gr{n}{\rd X}$, set of 1\nbd cells is $\gr{n+1}{\rd X}$, and set of 2\nbd cells is the quotient $\slice{\gr{n+2}{\rd{X}}}{\simeq}$.
A consequence of this fact is that, to prove some equivalences of round diagrams that only rely on properties of pasting and units in codimension 1 and 2, we can rely on established facts about bicategories, in particular, the celebrated coherence theorem \cite[Theorem 3.1]{maclane1963natural}, as well as the soundness of the calculus of string diagrams \cite{hinze2023introducing}.

The results of this section are not needed at any other point in the article; we include them to reassure the reader who may be wondering whether it is sound to use string-diagrammatic reasoning to prove facts about diagrammatic sets.
We note that \cite{rice2020coinductive}, for instance, treats soundness of planar isotopy of string diagrams as an axiom for any reasonable theory of higher categories.

\begin{dfn}[The bicategory of round $n$-dimensional diagrams]
	Let $X$ be a diagrammatic set, $n \in \mathbb{N}$.
	The \emph{bicategory of round $n$\nbd dimensional diagrams in $X$} is the bicategory $\bic{n}{X}$ determined by the following data.
	\begin{itemize}
		\item The set of 0\nbd cells is $\gr{0}{\bic{n}{X}} \eqdef \gr{n}{\rd X}$.
		\item The set of 1\nbd cells is $\gr{1}{\bic{n}{X}} \eqdef \gr{n+1}{\rd X}$, and the type of a 1\nbd cell $a$ is $\bd{}{-}a \celto \bd{}{+}a$.
		\item The set of 2\nbd cells is $\gr{2}{\bic{n}{X}} \eqdef \slice{\gr{n+2}{\rd X}}{\simeq}$, and the type of a 2\nbd cell $\isocl{t}$ is $\bd{}{-}t \celto \bd{}{+}t$; note that this is independent of the representative $t$.
		\item The horizontal composition of two 1\nbd cells $a\colon u \celto v$ and $b\colon v \celto w$ is $b * a \eqdef a \cp{} b\colon u \celto w$.
		\item The identity on a 0\nbd cell $u$ is $1_u \eqdef \un u\colon u \celto u$.
		\item The horizontal composition of two 2\nbd cells $\isocl{t}\colon a \celto b$ and $\isocl{s}\colon c \celto d$ is $\isocl{s} * \isocl{t} \eqdef \isocl{(t \cp{n} s) \cp{} \un{(b \cp{} d)}}\colon c * a \celto d * b$.
		\item The vertical composition of two 2\nbd cells $\isocl{t}\colon a \celto b$ and $\isocl{s}\colon b \celto c$ is $\isocl{s} \after \isocl{t} \eqdef \isocl{t \cp{} s}\colon a \celto c$.
		\item The identity on a 1\nbd cell $a\colon u \celto v$ is $1_a \eqdef \isocl{\un a}\colon a \celto a$.
		\item The associator indexed by three 1\nbd cells $a, b, c$ is an identity.
		\item The right unitor indexed by a 1\nbd cell $a\colon u \celto v$ is $r_a \eqdef \isocl{\lun{}a}\colon a \celto a * 1_u$.
		\item The left unitor indexed by a 1\nbd cell $a\colon u \celto v$ is $\ell_a \eqdef \isocl{\run{}a}\colon 1_v * a \celto a$.
	\end{itemize}
\end{dfn}

\begin{rmk}
	We use the classical order of composition in $\bic{n}{X}$, as opposed to the diagrammatic order of composition used for pasting, as an extra step to avoid confusion between the two.
	Note that this flips the side of unitors.
	We also make an arbitrary convenient choice for the default direction of unitors, which is not standard in the literature anyway.
\end{rmk}

\begin{prop} \label{prop:bicategory_well_defined}
	Let $X$ be a diagrammatic set.
	Then $\bic{n}X$ is well-defined as a bicategory.
\end{prop}
\begin{proof}
	To begin, vertical composition of 2\nbd cells is associative on the nose.
	Let $\isocl{t}\colon a \celto b$ be a 2\nbd cell; then the equivalences $t \cp{} \un{b} \simeq t \simeq \un{a} \cp{} t$ exhibit $1_b * \isocl{t} = \isocl{t} = \isocl{t} * 1_a$.
	Next, we prove naturality of horizontal composition with respect to vertical composition and units.
	Given 2\nbd cells $\isocl{t}\colon a \celto b$, $\isocl{s}\colon b \celto c$, $\isocl{t'}\colon a' \celto b'$, and $\isocl{s'}\colon b' \celto c'$ such that $\bd{}{+}a = \bd{}{-}a'$, we have
	\begin{align*}
		((t \cp{} s) \cp{n} (t' \cp{} s')) \cp{} \un{(c \cp{} c')} 
			& = (t \cp{n} t') \cp{} (s \cp{} s') \cp{} \un{(c \cp{} c')} \\
			& \simeq (t \cp{n} t') \cp{} \un{(b \cp{} b')} \cp{} (s \cp{n} s') \cp{} \un{(c \cp{} c')},
	\end{align*}
	exhibiting $(\isocl{s'} \after \isocl{t'}) * (\isocl{s} \after \isocl{t}) = (\isocl{s'} * \isocl{s}) \after (\isocl{t'} * \isocl{t})$.
	Given 1\nbd cells $a\colon u \celto v$ and $b\colon v \celto w$, we have
	\begin{align*}
		(\un{(a)} \cp{n} \un{(b)}) \cp{} \un{(a \cp{} b)} 
			& = \un{(a)} \cpsub{} (\un{(b)} \cpsub{} \un{(a \cp{} b)}) \\
			& \simeq \un{(b)} \cpsub{} \un{(a \cp{} b)} \\
			& \simeq \un{(a \cp{} b)},
	\end{align*}
	exhibiting $1_b * 1_a = 1_{b * a}$.
	Next, we prove that associators and unitors are natural in their parameters.
	This is automatic for associators since they are strict.
	Let $\isocl{t}\colon a \celto b$ be a 2\nbd cell, where $a, b\colon u \celto v$.
	Then
	\begin{align*}
		t \cp{} \lun{}b 
			& \simeq \lun{}a \subcp{} t \\
			& \simeq (\lun{}a \subcp{} t) \subcp{} \un{(\un{u})} \\
			& = \lun{}a \cp{} (\un{(\un{u})} \cp{n} t) \\
			& \simeq \lun{}a \cp{} (\un{(\un{u})} \cp{n} t) \cp{} \un{(\un{u} \cp{} b)},
	\end{align*}
	exhibiting $r_b \after t = (t * 1_{1_u}) \after r_a$, which proves that right unitors are natural in their parameter.
	The proof that left unitors are natural in their parameter is dual.
	Moreover, because units and unitors are weakly invertible in $X$, associators and unitors are componentwise invertible in $\bic{n}{X}$.
	Finally, the pentagon equation holds automatically for strict associators, whereas, given 1\nbd cells $a\colon u \celto v$ and $b\colon v \celto w$,
	\begin{align*}
		\un{(a \cp{} b)} 
			& = (a \cp{} \lun{}b)\cp{}(\run{}a \cp{} b) \\
			& \simeq \un{a} \cpsub{} (a \cp{} \lun{} b) \cp{} (\run{}a \cp{} b) \subcp{} \un{b} \\
			& = (\un{a} \cp{n} \lun{} b) \cp{} (\run{}a \cp{n} \un{b}) \\
			& \simeq (\un{a} \cp{n} \lun{} b) \cp{} (\run{}a \cp{n} \un{b}) \cp{} \un{(a \cp{} b)} \\
			& \simeq (\un{a} \cp{n} \lun{} b) \cp{} \un{(a \cp{} \un{v} \cp{} b)}
			\cp{} (\run{}a \cp{n} \un{b}) \cp{} \un{(a \cp{} b)}
	\end{align*}
	exhibits $1_{b * a} = (1_b * \ell_a) \after (r_b * 1_a)$, which is equivalent to the triangle equation in a strictly associative bicategory.
	This completes the proof.
\end{proof}

\begin{lem} \label{lem:equivalence_in_bicategory}
	Let $X$ be a diagrammatic set, $n \in \mathbb{N}$, and let $e \in \gr{n+1}{\rd X}$ and $h \in \gr{n+2}{\rd X}$.
	Then
	\begin{enumerate}
		\item $h$ is an equivalence in $X$ if and only if $\isocl{h}$ is invertible in $\bic{n}{X}$,
		\item $e$ is an equivalence in $X$ if and only if $e$ is an equivalence in $\bic{n}{X}$.
	\end{enumerate}
\end{lem}
\begin{proof}
	Straightforward using Corollary \ref{cor:equivalence_is_invertible_up_to_equivalence}.
\end{proof}

\noindent The following proof is an example of how one can leverage Proposition \ref{prop:bicategory_well_defined} to import known facts about bicategories into the theory of diagrammatic sets.

\begin{prop} \label{prop:adjointification}
	Let $X$ be a diagrammatic set, let $e \in \eqv X$ of type $u \celto v$, and let $e^*$ be a weak inverse of $e$.
	Then there exist invertors $z\colon e \cp{} e^* \celto \un u$ and $h\colon \un v \celto e^* \cp{} e$ that are ``adjoint up to equivalence'', that is, satisfy
	\[	(e \cp{} h) \cp{} (z \cp{} e) \simeq \run{}e \cp{} \lun{}e,
		\quad \quad
		(h \cp{} e^*) \cp{} (e^* \cp{} z) \simeq \rev{(\lun{}e^*)} \cp{} \rev{(\run{}e^*)}.
	\]
\end{prop}
\begin{proof}
	Let $n \eqdef \dim{u} = \dim{v}$.
	By Lemma \ref{lem:equivalence_in_bicategory}, $e$ is an equivalence in the bicategory $\bic{n}{X}$.
	By a standard result in bicategory theory \cite[Proposition 6.2.4]{johnson2021twodimensional}, $e$ is part of an adjoint equivalence exhibited by invertible 2\nbd cells $\isocl{z}\colon e^* \after e \celto 1_u$ and $\isocl{h}\colon 1_v \celto e \after e^*$ in $\bic{n}{X}$; that is, $\isocl{z}$ and $\isocl{h}$ satisfy
	\[
		(1_e * \isocl{z}) \after (\isocl{h} * 1_e) = r_e \after \ell_e, 
		\quad \quad
		(\isocl{z} * 1_{e^*}) \after (1_{e^*} * \isocl{h}) = \invrs{\ell_{e^*}} \after \invrs{r_{e^*}}.
	\]
	Translating to $X$ according to the definition, the first equation becomes
	\[
		(\un{e} \cp{n} h) \cp{} \un{(e \cp{} e^* \cp{} e)} \cp{} (z \cp{n} \un{e}) \cp{} \un{(\un{u} \cp{} e)} \simeq \run{}e \cp{} \lun{}e,
	\]
	whose left-hand side is equal to
	\begin{align*}
		\un{e} \cpsub{} (e \cp{} h) \cp{} &\un{(e \cp{} e^* \cp{} e)} \cp{} (z \cp{} e) \subcp{} \un{e} \cp{} \un{(\un{u} \cp{} e)} \\
					    & \simeq \un{e} \cpsub{} (e \cp{} h) \cp{} (z \cp{} e) \subcp{} \un{e} \cp{} \un{(\un{u} \cp{} e)} \\
					    & \simeq (e \cp{} h) \cp{} (z \cp{} e)
	\end{align*}
	using appropriate left and right unitors.
	The other equation is dual.
\end{proof}

%% file: division.tex
\section{The division lemma} \label{sec:division}

\subsection{Rounded higher contexts}

\begin{dfn}[Rounded higher contexts]
	Let $X$ be a diagrammatic set, $v, w$ be parallel round diagrams in $X$, and $\fun{F}$ be a round context on $\pd X(v, w)$.
	For each parallel pair $a, b\colon c \celto d$ in $\rd X(v, w)$, we define, inductively on dimension, a round context $\R_{a,b}\fun{F}$ on $\pd X(a, b)$.
	We let $\R_{v, w}\fun{F} \eqdef \fun{F}$, and
	\[
		\R_{a, b}\fun{F} \eqdef {(\R_{c,d}\fun{F})_{a,b}-} \cp{} \un (\R_{c,d}\fun{F}b)\colon 
		\pd X(a, b) \to \pd X(\R_{c,d}\fun{F}a, \R_{c,d}\fun{F}b).
	\]
	We call these the \emph{rounded higher contexts} associated with $\fun{F}$.
\end{dfn}

\begin{rmk}
	If $\fun{F}$ is weakly invertible, then each rounded higher context $\R_{a, b}\fun{F}$ is weakly invertible.
\end{rmk}

\begin{rmk}
	If $a, b$ are of type $c \celto d$, then $\R_{a,b}\fun{F} = \R_{a,b}(\R_{c,d}\fun{F})$.
\end{rmk}

\begin{lem} \label{lem:one_step_rounding_compsite_context}
	Let \( X \) be a diagrammatic set, \( \fun{F} \) a round context on $\pd X(v, w)$, \( a, b \colon v \celto w \) a parallel pair of round diagrams, and $\fun{G}$ a round context on $\pd X(\fun{F}a, \fun{F}b)$.
	If $\fun{G}\fun{F}_{a,b}$ is round, then $\fun{G}\fun{F}_{a,b}\simeq \fun{G}\R_{a,b}\fun{F}$.
\end{lem}
\begin{proof}
	Let $k \eqdef \dim{\fun{G}}$.
	By Lemma \ref{lem:context_factorisation}, we may write $\fun{G} = \fun{T}\fun{G'}_{\fun{F}a,\fun{F}b}$ where $\fun{T}$ is trim and round, and $\fun{G'}$ is round of dimension $k-1$.
	Moreover, by Lemma \ref{lem:trim_context_layering}, we may write $\fun{T} = \ell \cp{k-1} - \cp{k-1} r$ for a pair of pasting diagrams with $\dim{\ell}, \dim{r} \leq k$.
	Then we have the following sequence of natural equivalences of round contexts:
	\begin{align*}
		\fun{G}\fun{F}_{a,b} & = \ell \cp{k-1} {\fun{G'F}_{a,b}-} \cp{k-1} r \\
				     & \simeq \un{(\bd{k-1}{-}\ell)} \cp{} (\ell \cp{k-1} {\fun{G'F}_{a,b}-} \cp{k-1} r) 
				     	&& \text{by left unitor} \\
					& \simeq (\un{(\bd{k-1}{-}\ell)} \cp{k-1} \ell \cp{k-1} {\fun{G'F}_{a,b}-}) \subcp{} \un{(\fun{F}b)} \cp{k-1} r 
					&& \text{by right unitor} \\
					& \simeq \ell \cp{k-1} ({\fun{G'F}_{a,b}-} \subcp{} \un{(\fun{F}b)}) \cp{k-1} r 
					&& \text{by left unitor} \\
					& = \ell \cp{k-1} \fun{G'}_{\fun{F}a, \fun{F}b}({\fun{F}_{a,b}-} \cp{} \un{(\fun{F}b)}) \cp{k-1} r = \fun{G}\R_{a,b}\fun{F}.
	\end{align*}
	We conclude by Proposition \ref{prop:natural_equivalence_is_equivalence_relation}.
\end{proof}

\begin{lem} \label{lem:cylinder_eversion}
	Let $X$ be a diagrammatic set, $\fun{F}, \fun{G}$ be round contexts on $\pd X(v, w)$, and $\th\colon \fun{F} \eqvto \fun{G}$ be a natural equivalence.
	Then, for each parallel pair $a, b\colon c \celto d$ in $\rd X(v, w)$, there exist
	\begin{enumerate}
		\item a weakly invertible round context 
			\[ \fun{C}_{a,b}\colon \pd X(\R_{c,d}\fun{G}a, \R_{c,d}\fun{G}b) \to \pd X(\R_{c,d}\fun{F}a, \R_{c,d}\fun{F}b), \]
		\item a natural equivalence 
			\[ \th_{a,b}\colon \R_{a,b}\fun{F} \eqvto \fun{C}_{a,b}\R_{a,b}\fun{G} \]
		of round contexts on $\pd X(a, b)$,
	\end{enumerate}
	such that, letting $\fun{C}_{v, w} \eqdef -$ and $\th_{v, w} \eqdef \th$, we have, inductively,
	\[
		\fun{C}_{a, b} = \th_{c, d}a \cp{} 
			{\fun{C}_{c, d}-} \cp{} \th^*_{c, d}b.
	\]
\end{lem}
\begin{proof}
	We have, by assumption, $\th_{v, w}\colon \R_{v,w}\fun{F} \eqvto \R_{v,w}\fun{G}$.
	We may then assume, inductively, that we have defined $\th_{c,d}\colon \R_{c,d}\fun{F} \eqvto \fun{C}_{c,d}\R_{c,d}\fun{G}$.
	Then any choice of a componentwise weak inverse $\th^*_{c,d}$ allows us to define $\fun{C}_{a,b}$ with the specified type.
	We have the following sequence of natural equivalences, where we omit explicit promotions for the sake of readability:
	\begin{align*}
		\fun{C}_{a,b}\R_{a,b}\fun{G} 
			& = \th_{c, d}a \cp{} \fun{C}_{c, d}{\R_{a,b}(\R_{c,d}\fun{G})-} \cp{} \th^*_{c, d}b \\
			& \simeq \th_{c, d}a \cp{} \fun{C}_{c, d}{\R_{c,d} \fun{G}-} \cp{} \th^*_{c, d}b 
				&& \text{by Lemma \ref{lem:one_step_rounding_compsite_context}} \\
			& \simeq {\R_{c,d}\fun{F}-} \cp{} \th_{c, d}b \cp{} \th^*_{c, d} b  
				&& \text{by naturality of $\th_{c, d}$} \\
			& \simeq {\R_{c,d}\fun{F}-} \cp{} \un(\R_{c, d}\fun{F}b) 
				&& \text{by pfw with left invertor} \\
			& = \R_{a,b}\fun{F}.
	\end{align*}
	By Proposition \ref{prop:natural_equivalence_is_equivalence_relation}, this defines \( \th_{a,b} \colon \R_{a,b} \fun{F} \eqvto \fun{C}_{a,b}\R_{a,b}\fun{G} \).
\end{proof}

\begin{lem} \label{lem:rounding_unitality}
	Let $X$ be a diagrammatic set, let $v, w$ be parallel round diagrams in $X$, and let $\fun{I}$ denote the identity context on $\pd X(v, w)$.
	Then, for each parallel pair $a, b\colon c \celto d$ in $\rd X(v, w)$, there exist
	\begin{enumerate}
		\item a weakly invertible round context 
			\[ \fun{J}_{a,b}\colon \pd X(\R_{c,d}\fun{I}a, \R_{c,d}\fun{I}b) \to \pd X(a, b), \]
		\item a natural equivalence 
			\[ \eta_{a,b}\colon - \eqvto \fun{J}_{a,b}\R_{a,b}\fun{I} \]
		of round contexts on $\pd X(a, b)$,
	\end{enumerate}
	such that, letting $\fun{J}_{v, w} \eqdef -$ and $\eta_{v, w} \eqdef \un{}$, we have, inductively,
	\[
		\fun{J}_{a, b} = \eta_{c, d}a \cp{} 
			{\fun{J}_{c, d}-} \cp{} \eta^*_{c, d}b.
	\]
\end{lem}
\begin{proof}
	A straightforward variation on the proof of Lemma \ref{lem:cylinder_eversion}.
\end{proof}

\begin{lem} \label{lem:rounding_functoriality}
	Let $X$ be a diagrammatic set, $\fun{F}\colon \pd X(v, w) \to \pd X(v', w')$ a round context, and $\fun{G}$ a round context on $\pd X(v', w')$.
	Then, for each parallel pair $a, b\colon c \celto d$ in $\rd X(v, w)$, there exist
	\begin{enumerate}
		\item a weakly invertible round context $\fun{M}_{a,b}$ of type
			\[
				\pd X((\R_{c',d'}\fun{G})\R_{c,d}\fun{F}a, (\R_{c',d'}\fun{G})\R_{c,d}\fun{F}b))
				\to
				\pd X(\R_{c,d}(\fun{GF})a, \R_{c,d}(\fun{GF})b),
			\]
		\item a natural equivalence
			\[
				\mu_{a,b} \colon \R_{a,b}(\fun{GF}) \eqvto \fun{M}_{a,b}(\R_{a',b'}\fun{G})\R_{a,b}\fun{F}
			\]
			of round contexts on $\pd X(a, b)$,
	\end{enumerate}
	such that, letting $\fun{M}_{v,w} \eqdef -$ and $\mu_{v,w} \eqdef \un{}$, we have, inductively,
	\[
		\fun{M}_{a,b} = \mu_{c,d}a \cp{} {\fun{M}_{c,d}-} \cp{} \mu_{c,d}^* b, 
		\quad \quad
		a' = \R_{c,d}\fun{F}a, 
		\quad \quad b' = \R_{c,d}\fun{F}b.
	\]
\end{lem}
\begin{proof}
	Another easy variation on the proof of Lemma \ref{lem:cylinder_eversion}.
\end{proof}

\begin{rmk}
	The results of this section appear to be at least superficially related to \cite[Construction 3.2.2]{fujii2024omega}, with our weakly invertible contexts playing the role of ``paddings'', and the various natural equivalences establishing their naturality as in \cite[Lemma 3.2.8]{fujii2024omega}.
\end{rmk}

%%%%%%%%%%%%%

\subsection{Proof of the division lemma}

\begin{dfn}[Weakly unique solutions to equations]
	Let $X$ be a diagrammatic set, $\fun{F}\colon \pd X(v, w) \to \pd X(v', w')$ a round context, $b\colon v' \celto w'$ a round diagram, and $n \eqdef \dim\fun{F}$.
	A \emph{solution} to the equation $\fun{F}x \qeq b$ in the indeterminate $x \in \gr{n}{\rd X(v, w)}$ is a round diagram $a\colon v \celto w$ such that $\fun{F}a \simeq b$.
	A solution is \emph{weakly unique} if, for all parallel pairs of round diagrams $a, a'\colon v \celto w$, $\fun{F}a \simeq \fun{F}a'$ implies $a \simeq a'$.
\end{dfn}

\noindent The \emph{division lemma} is the following statement.

\begin{lem} \label{lem:division_lemma}
	Let $X$ be a diagrammatic set, $\fun{E}\colon \pd X(v, w) \to \pd X(v', w')$ a weakly invertible round context, $b\colon v' \celto w'$ a round diagram, and $n \eqdef \dim\fun{E}$.
	Then the equation $\fun{E}x \qeq b$ in the indeterminate $x \in \gr{n}{\rd X(v, w)}$ has a weakly unique solution.
\end{lem}

\begin{comm}
	The division lemma can be read as the statement that a weakly invertible round context $\fun{E}\colon \pd X(v, w) \to \pd X(v', w')$ establishes a bijection between $\gr{n}{\rd X(v, w)}$ and $\gr{n}{\rd X(v', w')}$ up to $(n+1)$\nbd dimensional equivalences.
\end{comm}

\begin{rmk}
	When $\fun{E} = {e \cpsub{\iota} -}$ or ${- \subcp{\iota} e}$ for some $e \in \gr{n}{\eqv X}$, then a solution to $\fun{E}x \qeq b$ exists by definition of equivalence.
\end{rmk}

\noindent The rest of the article will be devoted to the proof of the division lemma.

\begin{dfn}[Factorisation preorder on round contexts]
	Let $X$ be a diagrammatic set and $\fun{F}, \fun{G}$ be round contexts on $\pd X(v, w)$.
	We write $\fun{F} \fact \fun{G}$ if and only if there exists a round context $\fun{C}$ such that $\fun{G} \simeq \fun{CF}$.
	This determines a preorder, the \emph{factorisation preorder} on round contexts on $\pd X(v, w)$.
\end{dfn}

\begin{rmk}
	Given any round context $\fun{F}$, it is always the case that $- \fact \fun{F}$.
\end{rmk}

\begin{lem} \label{lem:factorisation_implies_kernel_containment}
	Let $X$ be a diagrammatic set, let $\fun{H}, \fun{K}$ be parallel round contexts with codomain $\pd X(v, w)$, let $\fun{F}, \fun{G}$ be round contexts on $\pd X(v, w)$, and suppose $\fun{F} \fact \fun{G}$.
	Then $\fun{FH} \simeq \fun{FK}$ implies $\fun{GH} \simeq \fun{GK}$.
\end{lem}
\begin{proof}
	Let $\fun{C}$ be a round context such that $\fun{G} \simeq \fun{CF}$.
	By closure of natural equivalences under left context, $\fun{FH} \simeq \fun{FK}$ implies $\fun{CFH} \simeq \fun{CFK}$, so by closure under right context $\fun{GH} \simeq \fun{GK}$.
\end{proof}

\begin{dfn}[Weak inverse of a round context]
	Let $X$ be a diagrammatic set and $\fun{E}\colon \pd X(v, w) \to \pd X(v', w')$ a round context.
	A \emph{weak inverse} of $\fun{E}$ is a round context $\fun{E}^*\colon \pd X(v', w') \to \pd X(v, w)$ such that $\fun{E}^*\fun{E} \simeq -$ and $\fun{E}\fun{E}^* \simeq -$.
\end{dfn}

\begin{lem} \label{lem:weakly_invertible_trim_context_has_weak_inverse}
	Let $X$ be a diagrammatic set, let $v, w$ be a parallel pair of round diagrams in $X$, and let $\fun{E}$ be a weakly invertible trim context on $\pd X(v, w)$.
	Then $\fun{E}$ has a weakly invertible weak inverse $\fun{E}^*$.
\end{lem}
\begin{proof}
	We proceed by structural induction on $\fun{E}$.
	If $\fun{E}$ is of the form $e \cpsub{\iota} -$ for some equivalence $e$ and $\iota\colon \bd{}{+}e \submol v$, let $e^*$ be a weak inverse of $e$, and let $\fun{E}^* \eqdef e^* \cpsub{j} -$, where $j$ is the inclusion of $\bd{}{+}e^* = \bd{}{-}e$ into $\subs{v}{\bd{}{-}e}{\iota(\bd{}{+}e)}$.
	Then
	\begin{align*}
		\fun{E}^*\fun{E} & = (e^* \cp{} e) \cpsub{\iota} - \\
			& \simeq \un{(\bd{}{+}e)} \cpsub{\iota} - 
				&& \text{by pushforward with invertor} \\
			& \simeq -
				&& \text{by left unitor},
	\end{align*}
	and similarly $\fun{E}\fun{E}^* \simeq -$.
	The proof in the case that $\fun{E}$ is of the form $- \subcp{\iota} e$ is dual.
	If $\fun{E} = -$ is the identity context, then $\fun{E}^* \eqdef -$ is a weak inverse of \( \fun{E} \) by Proposition \ref{prop:natural_equivalence_is_equivalence_relation}.
	Finally, if $\fun{E} = \fun{GF}$ with $\fun{F}, \fun{G}$ weakly invertible trim contexts, by the inductive hypothesis $\fun{F}$ and $\fun{G}$ have weakly invertible weak inverses $\fun{F}^*$, $\fun{G}^*$, respectively.
	Then, letting $\fun{E}^* \eqdef \fun{F}^*\fun{G}^*$, since $\simeq$ is a congruence, we have
	\[
		\fun{E}^*\fun{E} \simeq \fun{F}^*\fun{G}^*\fun{GF} \simeq \fun{F}^*\fun{F} \simeq -,
	\]
	and similarly $\fun{E}\fun{E}^* \simeq -$. 
	This concludes the proof.
\end{proof}

\begin{lem} \label{lem:if_factors_through_identity_so_do_roundings}
	Let $X$ be a diagrammatic set, let $\fun{F}$ be a round context on $\pd X(v, w)$, and suppose $\fun{F} \fact -$.
	Then, for each parallel pair $a,b\colon c \celto d$ in $\rd X(v, w)$, we have $\R_{a,b}\fun{F} \fact -$.
\end{lem}
\begin{proof}
	Let $\fun{I}$ be the identity context on $\pd X(v, w)$, and let $\fun{C}$ be a round context such that $\fun{I} \simeq \fun{CF}$.
	By Lemma \ref{lem:rounding_functoriality}, we have $\R_{a,b}\fun{F} \fact \R_{a,b}(\fun{CF})$.
	By Lemma \ref{lem:cylinder_eversion}, we have $\R_{a,b}(\fun{CF}) \fact \R_{a,b}\fun{I}$.
	Finally, by Lemma \ref{lem:rounding_unitality}, we have $\R_{a,b}\fun{I} \fact -$.
	We conclude by transitivity of the preorder.
\end{proof}

\begin{lem} \label{lem:weakly_invertible_factors_through_identity}
	Let $X$ be a diagrammatic set and let $\fun{E}$ be a weakly invertible round context on $\pd X(v, w)$.
	Then $\fun{E} \fact -$.
\end{lem}
\begin{proof}
	We proceed by induction on $k \eqdef \dim{\fun{E}}$.
	If $k = 1$, then $\fun{E}$ is trim, so it has a weak inverse $\fun{E}^*$ by Lemma \ref{lem:weakly_invertible_trim_context_has_weak_inverse}.
	Then $\fun{E}^*\fun{E} \simeq -$ exhibits $\fun{E} \fact -$.
	Suppose $\dim{\fun{E}} > 1$.
	By Lemma \ref{lem:context_factorisation}, we can write $\fun{E} = \fun{T}\fun{F}_{v, w}$ where $\fun{T}, \fun{F}$ are weakly invertible, $\fun{T}$ is trim, and $\fun{F}$ is round with $\dim{\fun{F}} = k - 1$.
	By Lemma \ref{lem:one_step_rounding_compsite_context}, since $\fun{E}$ is round, we have $\fun{E} \simeq \fun{T}\R_{v,w}\fun{F}$.
	Let $\fun{T}^*$ be a weak inverse of $\fun{T}$; then $\fun{T}^*\fun{E} \simeq \R_{v, w}\fun{F}$ exhibits $\fun{E} \fact \R_{v,w}\fun{F}$.
	By the inductive hypothesis, $\fun{F} \fact -$, so by Lemma \ref{lem:if_factors_through_identity_so_do_roundings}, $\R_{v,w}\fun{F} \fact -$.
	We conclude by transitivity of the preorder.
\end{proof}

\begin{lem} \label{lem:left_inv_iff_right_inv_invertible_context}
	Let \( X \) be a diagrammatic set and let \( \fun{F}, \fun{G} \) be round contexts such that \( \fun{GF} \simeq - \). 
	If \( \fun{G} \) is weakly invertible, then it is a weak inverse of \( \fun{F} \).
\end{lem}
\begin{proof}
	From \( \fun{GF} \simeq - \), it follows that \( \fun{GFG} \simeq \fun{G} \).
	Since $\fun{G}$ is weakly invertible, by Lemma \ref{lem:weakly_invertible_factors_through_identity} we have $\fun{G} \fact -$, so by Lemma \ref{lem:factorisation_implies_kernel_containment} we conclude that $\fun{FG} \simeq -$.
\end{proof}

\begin{lem} \label{lem:roundings_preserve_weak_inverse}
	Let $X$ be a diagrammatic set, let $\fun{F}$ be a round context on $\pd X(v, w)$, and suppose that $\fun{F}$ has a weakly invertible weak inverse.
	Then, for each parallel pair $a,b\colon c \celto d$ in $\rd X(v, w)$, the context $\R_{a,b}\fun{F}$ has a weakly invertible weak inverse.
\end{lem}
\begin{proof}
	Let $\fun{I}$ denote the identity context on $\pd X(v, w)$ and let $\fun{G}$ be a weakly invertible round context such that $\fun{GF} \simeq \fun{I}$.
	Then, by Lemma \ref{lem:rounding_unitality}, Lemma \ref{lem:cylinder_eversion}, and Lemma \ref{lem:rounding_functoriality}, respectively, we have
	\[
		- \simeq \fun{J}_{a,b}\R_{a,b}\fun{I} \simeq \fun{J}_{a,b}\fun{C}_{a,b}\R_{a,b}(\fun{GF}) \simeq \fun{J}_{a,b}\fun{C}_{a,b}\fun{M}_{a,b}(\R_{a',b'}\fun{G})\R_{a,b}\fun{F}
	\]
	for some weakly invertible round contexts $\fun{J}_{a,b}$, $\fun{C}_{a,b}$, and $\fun{M}_{a,b}$.
	Then the composite $\fun{J}_{a,b}\fun{C}_{a,b}\fun{M}_{a,b}\R_{a',b'}\fun{G}$ is weakly invertible, so by Lemma \ref{lem:left_inv_iff_right_inv_invertible_context} it is a weak inverse of $\R_{a,b}\fun{F}$.
\end{proof}

\begin{thm} \label{thm:invertible_context_are_invertible}
	Every weakly invertible round context has a weakly invertible weak inverse.
\end{thm}
\begin{proof}
	Let $X$ be a diagrammatic set and let $\fun{E}$ be a weakly invertible round context on $\pd X(v, w)$.
	We proceed by induction on $k \eqdef \dim{\fun{E}}$.
	If $k = 1$, then $\fun{E}$ is trim, so by Lemma \ref{lem:weakly_invertible_trim_context_has_weak_inverse} it has a weakly invertible weak inverse.
	Suppose $k > 1$.
	By Lemma \ref{lem:context_factorisation}, we can write $\fun{E} = \fun{T}\fun{F}_{v, w}$ where $\fun{T}, \fun{F}$ are weakly invertible, $\fun{T}$ is trim, and $\fun{F}$ is round with $\dim{\fun{F}} = k - 1$.
	By Lemma \ref{lem:one_step_rounding_compsite_context}, since $\fun{E}$ is round, we have $\fun{E} \simeq \fun{T}\R_{v,w}\fun{F}$.
	Since $\fun{T}$ is trim, it has a weakly invertible weak inverse $\fun{T}^*$.
	By the inductive hypothesis, $\fun{F}$ has a weakly invertible weak inverse, so by Lemma \ref{lem:roundings_preserve_weak_inverse} $\R_{v,w}\fun{F}$ also has a weakly invertible weak inverse $\fun{G}$.
	Then $\fun{E}^* \eqdef \fun{G}\fun{T}^*$ is a weakly invertible weak inverse of $\fun{E}$.
\end{proof}

\begin{proof}[Proof of Lemma \ref{lem:division_lemma}]
	Let \( \fun{E}^* \) be a weak inverse of \( \fun{E} \).
	Then $\fun{E}^*b$ is a weakly unique solution to the equation $\fun{E}x \qeq b$.
\end{proof}